\documentclass[12pt,twoside]{article}
\usepackage{amsmath,amsfonts,amssymb,a4,enumerate,epsfig,array,theorem,psfrag}

\newtheorem{theo}{Theorem}

\newtheorem{prop}{Proposition}[section]
\newtheorem{coro}[prop]{Corollary}

\newtheorem{lemma}[prop]{Lemma}

\newcommand{\meio}{{\frac{1}{2}}}

\newcommand{\bfh}{{\mathbf h}}
\newcommand{\bfi}{{\mathbf i}}
\newcommand{\bfih}{\hat{\mathbf i}}
\newcommand{\bfj}{{\mathbf j}}
\newcommand{\bfk}{{\mathbf k}}
\newcommand{\bfkh}{\hat{\mathbf k}}

\newcommand{\bfn}{{\mathbf n}}
\newcommand{\bfs}{{\mathbf s}}
\newcommand{\bft}{{\mathbf t}}
\newcommand{\bfone}{{\mathbf 1}}

\newcommand{\Up}{\textit{Up}_3^{+}}
\newcommand{\Upum}{\textit{Up}_3^{1}}

\newcommand{\Next}{\operatorname{ns}}
\newcommand{\Free}{\operatorname{Free}}

\newcommand{\ZZ}{{\mathbb{Z}}}

\newcommand{\RR}{{\mathbb{R}}}

\newcommand{\Ss}{{\mathbb{S}}}

\newcommand{\DD}{{\mathbb{D}}}
\newcommand{\PP}{{\mathbb{P}}}
\newcommand{\HH}{{\mathbb{H}}}

\newcommand{\cH}{{\cal H}}
\newcommand{\cL}{{\cal L}}
\newcommand{\cI}{{\cal I}}
\newcommand{\cF}{{\cal F}}

\newcommand{\cC}{{\cal C}}
\newcommand{\cD}{{\cal D}}

\newcommand{\cM}{{\cal M}}
\newcommand{\cY}{{\cal Y}}

\newcommand{\cE}{{\cal E}}
\newcommand{\fF}{{\mathfrak F}}

\newcommand{\Bruhat}{\operatorname{Bru}}

\newcommand{\tot}{\operatorname{tot}}

\newcommand{\nobf}{\noindent\bf}

\def\qed{\unskip\nobreak\hfil\penalty50\hskip1.75em\null\nobreak\hfil
$\blacksquare$ {\parfillskip=0pt \finalhyphendemerits=0 \par}\goodbreak}

\begin{document}
\title{The homotopy type of spaces of \\
locally convex curves in the sphere}
\author{Nicolau C. Saldanha}
\maketitle

\begin{abstract}
A smooth curve $\gamma: [0,1] \to \Ss^2$ is
locally convex if its geodesic curvature is positive at every point.
J.~A.~Little showed that the space
of all locally convex curves $\gamma$
with $\gamma(0) = \gamma(1) = e_1$ and $\gamma'(0) = \gamma'(1) = e_2$
has three connected components
$\cL_{-\bfone,c}$, $\cL_{+\bfone}$, $\cL_{-\bfone,n}$.
The space $\cL_{-\bfone,c}$ is known to be contractible.
We prove that $\cL_{+\bfone}$ and $\cL_{-\bfone,n}$
are homotopy equivalent 
to $(\Omega\Ss^3) \vee \Ss^2 \vee \Ss^6 \vee \Ss^{10} \vee \cdots$
and $(\Omega\Ss^3) \vee \Ss^4 \vee \Ss^8 \vee \Ss^{12} \vee \cdots$,
respectively.
As a corollary, we deduce the homotopy type of the components
of the space $\Free(\Ss^1,\Ss^2)$ of free curves $\gamma: \Ss^1 \to \Ss^2$
(i.e., curves with nonzero geodesic curvature).
We also determine the homotopy type of the spaces
$\Free([0,1], \Ss^2)$ with fixed initial and final frames.
\end{abstract}

\section{Introduction}

\footnote{2000 {\em Mathematics Subject Classification}.
Primary 53C42; Secondary 34B05, 57N20, 53C29.
{\em Keywords and phrases} Convex curves,
topology in infinite dimension, 
periodic solutions of linear ODEs.}

\begin{figure}[ht]
\begin{center}
\epsfig{height=40mm,file=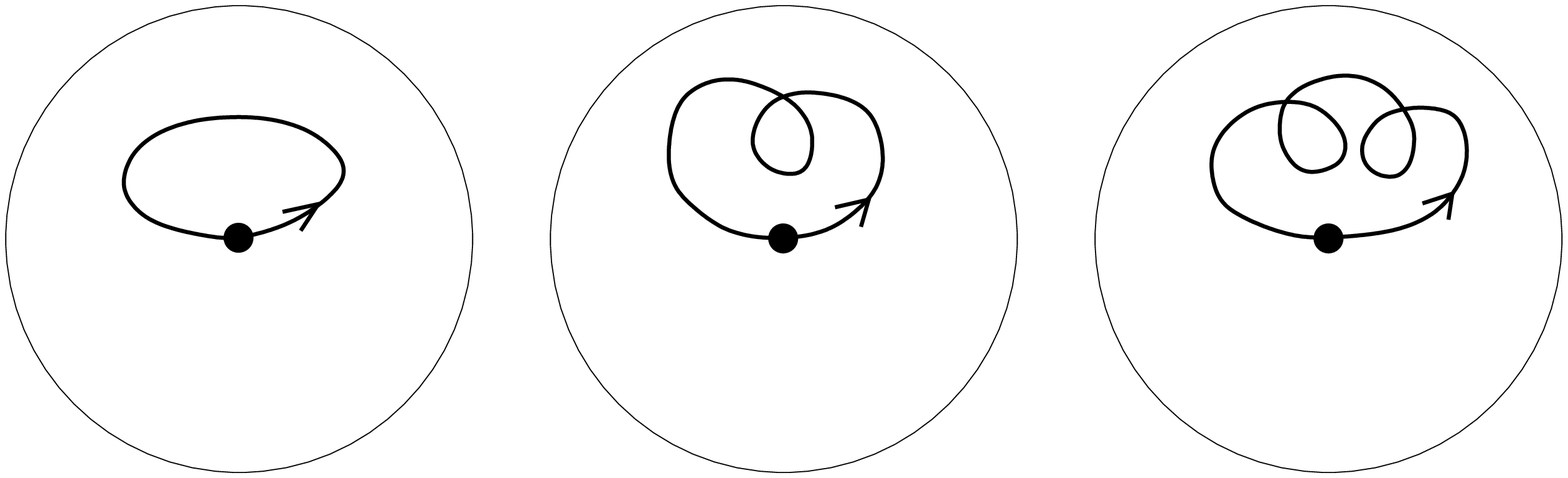}
\end{center}
\caption{Curves in $\cL_{-\bfone,c}$, $\cL_{+\bfone}$ and $\cL_{-\bfone,n}$.}
\label{fig:3comp}
\end{figure}

A curve $\gamma: [0,1] \to \Ss^2$ is called \textit{locally convex}
if its geodesic curvature is always positive,
or, equivalently, if $\det(\gamma(t), \gamma'(t), \gamma''(t)) > 0$ for all $t$.
Let $\cL_I$ be the space of all locally convex curves $\gamma$
with $\gamma(0) = \gamma(1) = e_1$ and $\gamma'(0) = \gamma'(1) = e_2$;
the precise topology for this space of curves
will be discussed in the paper.
J.~A.~Little \cite{Little} showed that $\cL_I$
has three connected components
$\cL_{-\bfone,c}$, $\cL_{+\bfone}$, $\cL_{-\bfone,n}$;
examples of curves in each connected component
are shown in Figure \ref{fig:3comp}.

The connected component $\cL_{-\bfone,c}$ can be defined to
be the set of simple curves in $\cL_I$:
the space $\cL_{-\bfone,c}$ is known to be contractible
(\cite{Anisov} and \cite{ShapiroM}, Lemma 5).
The aim of this paper is to determine the homotopy type
of the two remaining spaces $\cL_{+\bfone}$ and $\cL_{-\bfone,n}$.
Our main result is the following.

\begin{theo}
\label{theo:main}
The components $\cL_{+\bfone}$ and $\cL_{-\bfone,n}$
are homotopically equivalent to
$(\Omega\Ss^3) \vee \Ss^2 \vee \Ss^6 \vee \Ss^{10} \vee \cdots$
and $(\Omega\Ss^3) \vee \Ss^4 \vee \Ss^8 \vee \Ss^{12} \vee \cdots$,
respectively.
\end{theo}

Here $\Omega\Ss^3$ is the space of loops in $\Ss^3$,
i.e., the set of continuous maps $\alpha: [0,1] \to \Ss^3$
with $\alpha(0) = \alpha(1) = \bfone$, where $\bfone \in \Ss^3$ is a base point,
with the $C^0$ topology.
A more careful description of the connected components
$\cL_{+\bfone}$ and $\cL_{-\bfone,n}$ is given below.

A motivation for considering these spaces comes from differential equations.
Consider the linear ODE of order $3$:
\[ u'''(t) + h_1(t) u'(t) + h_0(t) u(t) = 0, \quad t \in [0,1]; \]
the set of pairs of potentials $(h_0,h_1)$ for which the equation
admits 3 linearly independent periodic solutions is homotopically
equivalent to $\cL_I$.
The corresponding problem in order $2$ is much simpler
(\cite{BST1}, \cite{BST2}, \cite{ST}).

\medskip

Alternatively, in Gromov's language (\cite{Gromov}, \cite{EM}),
given two smooth Riemannian manifolds $V^n$ and $W^q$,
a map $f: V^n \to W^q$ is \textit{free} (or \textit{second order free})
if the second order osculating space is non-degenerate
(we use covariant derivatives);
let $\Free(V,W)$ be the space of such free maps.
Perhaps the simplest non-trivial example here is $\Free(\Ss^1,\Ss^2)$,
the space of curves $\gamma: \Ss^1 \to \Ss^2$
with $\det(\gamma(t), \gamma'(t), \gamma''(t)) \ne 0$ for all $t$.
Thus $\Free(\Ss^1,\Ss^2)$ differs from our space $\cL_{I}$
only by the fact that in $\Free(\Ss^1,\Ss^2)$ negative determinants are allowed,
as are arbitrary initial frames $Q \in SO_3$.
The following result is a direct consequence
of Theorem \ref{theo:main}.

\begin{coro}
\label{coro:free}
The space $\Free(\Ss^1,\Ss^2)$ has six connected components,
with two homotopically equivalent to each of the following spaces:
\begin{gather*}
SO_3 \times \cL_{-\bfone,c} \approx
SO_3, \\
SO_3 \times \cL_{+\bfone} \approx
SO_3 \times \left((\Omega\Ss^3) \vee
\Ss^2 \vee \Ss^6 \vee \Ss^{10} \vee \cdots\right), \\
SO_3 \times \cL_{-\bfone,n} \approx
SO_3 \times \left((\Omega\Ss^3) \vee
\Ss^4 \vee \Ss^8 \vee \Ss^{12} \vee \cdots\right).
\end{gather*}
\end{coro}

Recall that if $q > \frac{n(n+3)}{2}$ then free maps
satisfy the parametrical h-principle;
we are here in the critical case
$n = 1$ and $q = \frac{n(n+3)}{2} = 2$ and 
the principle (incorrectly applied) would predict a (wrong) simpler answer.
This paper does not require any familiarity with these ideas
but the reader may notice that ideas similar to the h-principle
will play an important part.

\medskip

These spaces and variants have been discussed, among others,
by B.~Shapiro, M.~Shapiro and B.~Khesin (\cite{Shapiro2}, \cite{SK}).
These spaces are also the orbits of the second Gel'fand-Dikki brackets
and therefore have a natural symplectic structure (\cite{GD1}, \cite{GD2}).
Furthermore, these spaces are related to the orbit classification
of the Zamolodchikov Algebra (\cite{KS}, \cite{OK});
these interpretations shall not be used or discussed in this paper.

\medskip

Although the above authors ponder about the interest
of understanding the topology of such spaces,
their results deal mostly with $\pi_0$, i.e.,
with counting and identifying connected components.
The present author has also previously proved some weaker results
about the topology of these spaces,
regarding the fundamental group and the first few (co)homology groups
(\cite{S1}, \cite{S2});
these results are now of course
easy consequences of Theorem \ref{theo:main}.

\begin{coro}
\label{coro:oldtheo}
The spaces $\cL_{+\bfone}$ and $\cL_{-\bfone,n}$ are connected
and simply connected and, for $k > 0$,
their cohomology is given by
\[ H^k(\cL_{+\bfone};\ZZ) = \begin{cases}
0, & k \textrm{ odd}, \\
\ZZ^2, & 4 \mid (k + 2), \\
\ZZ, & 4 \mid k; \end{cases}
\quad
H^k(\cL_{-\bfone,n};\ZZ) = \begin{cases}
0, & k \textrm{ odd}, \\
\ZZ, & 4 \mid (k + 2), \\
\ZZ^2, & 4 \mid k. \end{cases} \]
\end{coro}




\bigskip

Let $\cI$ be the space
of immersions $\gamma: [0,1] \to \Ss^2$
(of class $C^k$ for some $k \ge 2$) with $\gamma'(t) \ne 0$,
$\gamma(0) = e_1$, $\gamma'(0) = e_2$.
Let $\cL \subset \cI$ be the subspace of locally convex curves;
thus, for $\gamma \in \cI$ we have $\gamma \in \cL$ 
if and only if $\det(\gamma(t), \gamma'(t), \gamma''(t)) > 0$ for all $t$.
For each $\gamma \in \cI$, consider its \textit{Frenet frame}
$\fF_\gamma: [0,1] \to SO_3$ defined by
\[ \begin{pmatrix} \gamma(t) & \gamma'(t) & \gamma''(t) \end{pmatrix} =
\fF_\gamma(t) R(t), \]
$R(t)$ being an upper triangular matrix with
$(R(t))_{11} > 0$ and $(R(t))_{22} > 0$
(the left hand side is the $3\times3$ matrix with columns
$\gamma(t)$, $\gamma'(t)$ and $\gamma''(t)$).
In other words, the first column of $\fF_\gamma(t)$ is $\gamma(t)$,
the second is the unit tangent vector
$\bft_\gamma(t) = \gamma'(t)/|\gamma'(t)|$
and the third column (which is now uniquely determined)
is the unit normal vector
$\bfn_\gamma(t) = \gamma(t) \times \bft_\gamma(t)$.
For $Q \in SO_3$, let $\cI_Q \subset \cI$ be the set of curves
$\gamma \in \cI$ for which $\fF_\gamma(1) = Q$;
similarly, let $\cL_Q = \cL \cap \cI_Q$.

The universal (double) cover of $SO_3$ is $\Ss^3 \subset \HH$,
the group of quaternions of absolute value $1$;
let $\Pi: \Ss^3 \to SO_3$ be the canonical projection.
For $\gamma \in \cI$,
the curve $\fF_\gamma$ can be lifted to define
$\tilde\fF_\gamma: [0, 1] \to \Ss^3$ with
$\tilde\fF_\gamma(0) = \bfone$, $\Pi \circ \tilde\fF_\gamma = \fF_\gamma$.
The value of $\tilde\fF_\gamma(1)$ partitions $\cI_Q$
as a disjoint union $\cI_{z} \sqcup \cI_{-z}$:
here $\Pi(\pm z) = Q$ and
$\gamma \in \cI_{z}$ if and only if $\tilde\fF_\gamma(1) = z$;
similarly, let $\cL_z = \cL_{\Pi(z)} \cap \cI_z$.
Notice that if $\gamma$ is a simple curve in $\cI_I$ then 
$\tilde\fF_\gamma(1) = -\bfone$ and therefore $\gamma \in \cI_{-\bfone}$.
We can now more precisely describe the three connected components of $\cL_I$:
the component $\cL_{+\bfone}$ is a special case of this definition
and $\cL_{-\bfone}$ has two connected components $\cL_{-\bfone,c}$
(convex or simple curves) and $\cL_{-\bfone,n}$ (non-simple).

A locally convex curve $\gamma: [t_0, t_1] \to \Ss^2$
is \textit{convex} if $\gamma$ intersects any geodesic
(great circle) at most twice. 
This definition requires a couple of clarifications.
First, endpoints do not count as intersections
so, for instance, a simple closed locally convex curve is convex.
Second, intersections are counted with multiplicity
so that a tangency counts as two intersections.
It follows from the definition that convex curves are simple.
In this paper we will see other equivalent definitions.

A matrix $Q \in SO_3$ is \textit{convex}
if there exists a convex arc $\gamma \in \cL$
with $\fF_\gamma(1) = Q$.
We shall see that the subset of $SO_3$ of convex matrices
is the disjoint union of one of the top dimensional Bruhat cells
and a few lower dimensional cells contained in its closure.
Similarly, a quaternion $z \in \Ss^3$ is \textit{convex}
if there exists a convex arc $\gamma \in \cL$
with $\tilde\fF_\gamma(1) = z$.
It is not hard to see that if $-z$ is convex then $z$ is not.

We can now state a more general version of our main theorem.
Here $\bfi \in \Ss^3$ is the usual quaternion.
Notice that $-\bfone$ is convex but $\bfone$, $\bfi$ and $-\bfi$ are not.

\begin{theo}
\label{theo:mainplus}
Let $z \in \Ss^3$.
Then the space $\cL_z$ is homotopically equivalent to
$\cL_{-\bfone}$ if $z$ is convex,
$\cL_{\bfone}$ if $-z$ is convex
and $\cL_{\bfi}$ otherwise.
Moreover, the following homotopy equivalences hold:
\[
\cL_{-\bfone} \approx
(\Omega\Ss^3) \vee \Ss^0 \vee \Ss^4 \vee \Ss^8 \vee \cdots; \quad
\cL_{+\bfone} \approx
(\Omega\Ss^3) \vee \Ss^2 \vee \Ss^6 \vee \Ss^{10} \vee \cdots; \quad
\cL_{\bfi} \approx \Omega\Ss^3.
\]
\end{theo}

We prove Theorem \ref{theo:mainplus} not just for the sake
of proving a stronger statement but mainly because
it is not clear how to produce a complete proof 
of Theorem \ref{theo:main}
(whose statement is simpler and more natural)
without strong use of Bruhat cells and other
algebraic notions.
By the time these ideas have been mastered,
both theorems are proved simultaneously.

Let $\Omega_z\Ss^3$ be the space of continuous curves
$\alpha: [0,1] \to \Ss^3$, $\alpha(0) = \bfone$, $\alpha(1) = z$: 
this is easily seen to be homeomorphic to $\Omega\Ss^3$
and the two spaces shall from now on be identified.
We just constructed the map $\tilde\fF: \cI_{z} \to \Omega_z\Ss^3$,
$\gamma \mapsto \tilde\fF_\gamma$.
It is a well-known fact
(which follows from the Hirsch-Smale Theorem)
that this map is a homotopy equivalence
(\cite{Morse}, \cite{Hirsch}, \cite{Smale}).
Theorem \ref{theo:mainplus} implies that
the inclusions $i: \cL_{z} \to \cI_{z}$
are homotopy equivalences only for certain quaternions $z$.

We now proceed to give an overview of the paper
and of the proof of the main theorems.
Section \ref{sect:topology} addresses the rather technical
issue of what, precisely, is the best topology for the space $\cL_{z}$.
As we shall see, we may allow for the juxtaposition of curves
(provided their Frenet frames agree);
on the other hand, when desirable,
we may assume curves to be smooth.
In Section \ref{sect:bruhat} we apply the construction
of Bruhat cells to our scenario.
We also study projective transformations.
The short Section \ref{sect:total}
collects a few useful facts
about the total (euclidean) curvature of spherical curves.

It follows from Little's results that
a circle drawn twice and a circle drawn four times
are in the same connected component of $\cL_I$.
In Section \ref{sect:2to4} we give a careful description
of a path joining these two curves.
We also prove a few facts about this path
which will be needed later.

As we have just mentioned,
the inclusion $i: \cL_{z} \to \cI_{z}$
need not be a homotopy equivalence.
In Section \ref{sect:homotosur} we see that
half of the story, so to speak, still holds.

\begin{prop}
\label{prop:homotosur}
Let $z \in \Ss^3$.
For any compact space $K$ and any function
$f: K \to \cI_{z}$ there exists $g: K \to \cL_{z}$
and a homotopy $H: [0,1] \times K \to \cI_{z}$
with $H(0,p) = f(p)$ and $H(1,p) = g(p)$ for all $p \in K$.
\end{prop}

The maps $i: \cL_{z} \to \cI_{z}$
therefore induce surjective maps
$\pi_k(\cL_{z}) \to \pi_k(\cI_{z})$.

\begin{figure}[ht]
\begin{center}
\epsfig{height=35mm,file=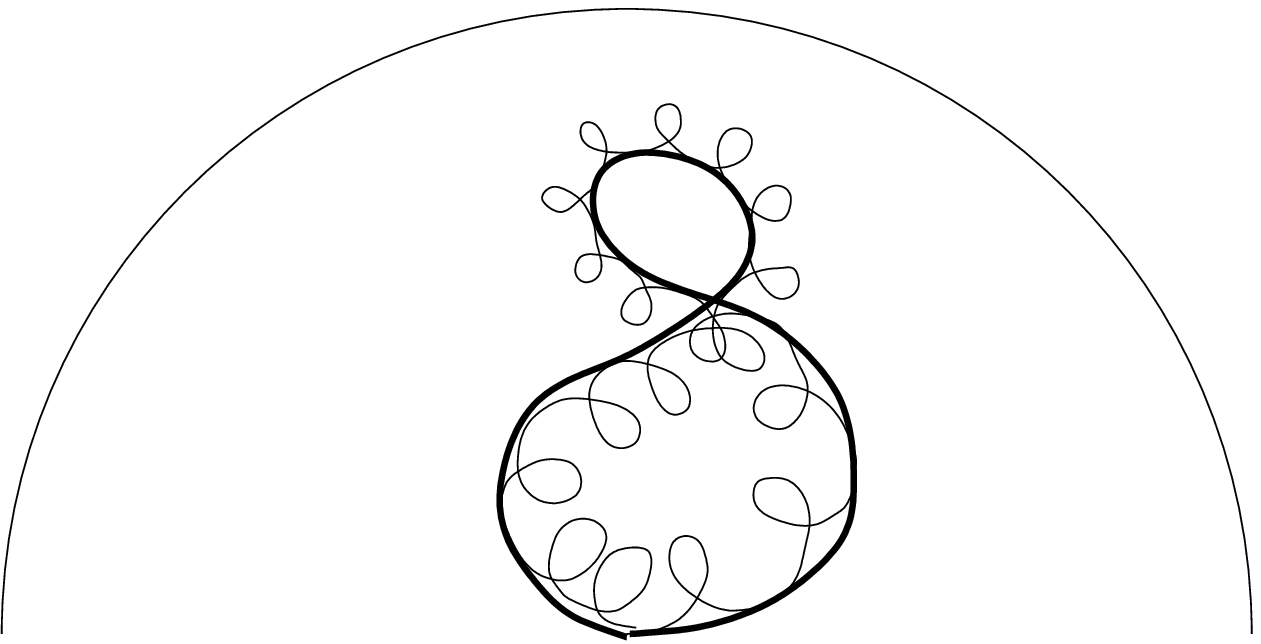}
\end{center}
\caption{Curves $\gamma_0 \in \cI_{\pm \bfone}$ (thick)
and $\gamma_1 \in \cL_{\pm \bfone}$ (thin).}
\label{fig:addloop}
\end{figure}

In a nutshell, if a curve has very large (positive)
geodesic curvature it looks like a phone wire and
we call it \textit{sufficiently loopy}.
Any compact family $\alpha_0: K \to \cI_z$ 
can be approximated in the $C^0$ topology by
a family $\alpha_1: K \to \cI_z$
of sufficiently loopy (and therefore locally convex) curves
(in Figure \ref{fig:addloop}, $\gamma_0 = \alpha_0(p)$
and $\gamma_1 = \alpha_1(p)$ for some $p \in K$).
Also, families of sufficiently loopy curves can be deformed
without losing the property of being locally convex.

The difficulty in proving (the false fact) that
the inclusion $\cL_z \subset \cI_z$ is a homotopy equivalence
is that there is no \textit{uniform} procedure
to add loops to locally convex curves
within the set of locally convex curves.

\begin{figure}[ht]
\centerline{\epsfig{height=35mm,file=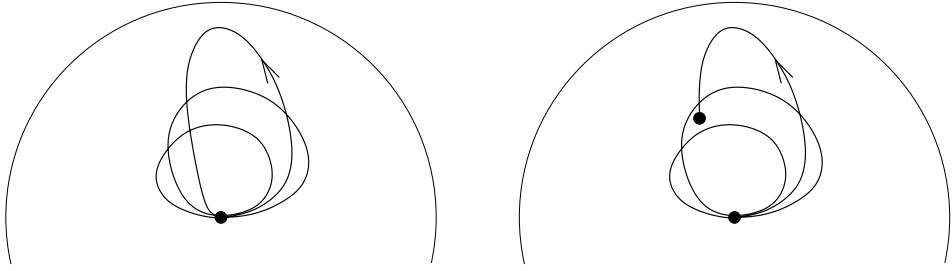}}
\caption{Two multiconvex curves of multiplicity $3$}
\label{fig:multiconvex}
\end{figure}

In Section \ref{sect:onion} we introduce a crucial construction
in our discussion:
a curve $\gamma \in \cL_{Q}$ is \textit{multiconvex}
of multiplicity $k$ if it is the juxtaposition
of $k-1$ simple closed convex curves with
a final $k$-th convex curve in $\cL_{Q}$
(see Figure \ref{fig:multiconvex}).

We prove in Lemma \ref{lemma:Mkcontract} that, given $Q$,
the closed subset $\cM_{k} \subset \cL_{Q}$
of multiconvex curves of multiplicity $k$ is
either empty or a contractible submanifold of codimension $2k-2$
with trivial normal bundle.
Assuming $-z$ convex,
we next construct in Lemma \ref{lemma:h2k-2}
maps $h_{2k-2}: \Ss^{2k-2} \to \cL_{(-\bfone)^k z}$ which
intersect $\cM_{k}$ transversally and exactly once
and are homotopic to a constant as maps $\Ss^{2k-2} \to \cI_{(-\bfone)^k z}$
(details of the construction of the path in Section \ref{sect:2to4}
are used here to verify that $h_{2k-2}$ has the desired properties). 
Intersection with $\cM_k$ shows that $h_{2k-2}$ is not homotopic
to a constant as a map $\Ss^{2k-2} \to \cL_I$
and therefore defines a non-trivial element of
the homotopy group $\pi_{2k-2}(\cL_{(-\bfone)^k z})$
which is taken to zero by the inclusion in $\cI_{(-\bfone)^k z}$.
Furthermore, intersection with $\cM_k$
defines an element of the cohomology group $H^{2k-2}(\cL_{(-\bfone)^k z})$
not in the image of $i^\ast: H^{\ast}(\cI_{(-\bfone)^k z}) \to
H^{\ast}(\cL_{(-\bfone)^k z})$
(compare with \cite{S1} and \cite{S2}).
This does not prove our main theorem yet
but already shows that if either $z$ or $-z$ is convex then
$\cL_z$ is not homotopically equivalent to $\cI_z$.

In Section \ref{sect:graft} we introduce \textit{grafting},
a process under which loops can sometimes be added to curves.
In Section \ref{sect:nextstep} we define the \textit{next step}
function and learn to tell apart \textit{good} and \textit{bad steps};
Bruhat cells are essential here.
Let $\cY_{z} = \cL_{z} \smallsetminus \bigcup_k \cM_{k}$
be the set of \textit{complicated} (i.e., not multiconvex) curves.
The new tools introduced
in Sections \ref{sect:graft} and \ref{sect:nextstep}
are then used in Section \ref{sect:Y} 
to understand the spaces $\cY_{z}$.

\begin{prop}
\label{prop:YI}
The inclusion $\cY_{z} \subset \cI_{z}$
is a weak homotopy equivalence.
\end{prop}



In order to determine the homotopy type of $\cL_z$,
start with $\cY_z$ and, for each $k$, add the set $\cM_k$.
It follows from what we have proved that adding $\cM_k$ 
(if nonempty) is equivalent to attaching a sphere $\Ss^{2k-2}$.
These properties are then sufficient to complete, in Section \ref{sect:theo},
the proof of Theorems \ref{theo:main} and \ref{theo:mainplus}.

The author would like to thank Dan Burghelea, Boris Khesin,
Boris Shapiro and Pedro Z\"uhlke for helpful conversations
and the referee and editor for valuable contributions;
thanks also go to the referee of \cite{S1} and \cite{S2}
for encouragement towards considering the problem discussed in this paper.
The author thanks the kind hospitality of
the Mathematics Department of The Ohio State University
during the winter quarters of 2004 and 2009 and of 
the Mathematics Department of the Stockholm University
during his visits to Sweden in 2005 and 2007.  
The author also gratefully acknowledges
the support of CNPq, Capes and Faperj (Brazil).

\section{Topology of $\cL_Q$}

\label{sect:topology}

We now attack a rather technical problem:
defining the best topological structure for the spaces
$\cL_Q$ and $\cI_Q$.

It turns out that different topological
structures obtain different spaces which are however
homotopically equivalent.
The $C^k$ metric for some $k \ge 2$ is a rather natural
choice but it has the inconvenience that when constructing
a homotopy we prefer not to be distracted by the necessity
to smoothen out certain points of our curves
(we want to be allowed, for instance, to consider
a curve $\gamma$ which is a juxtaposition of arcs of circle).

Given a smooth immersion $\gamma: [0,1] \to \Ss^2$,
let $\gamma'(t) = v_\gamma(t) \bft_\gamma(t)$
where $v_\gamma(t) = |\gamma'(t)| > 0$ and $\bft_\gamma(t)$
is the unit tangent vector to $\gamma$.
Let $\bfn_\gamma(t) = \gamma(t) \times \bft_\gamma(t)$
be the unit normal vector to $\gamma$, so that
\[ \bft_\gamma'(t) =
- v_\gamma(t) \gamma(t) + \hat v_\gamma(t) \bfn_\gamma(t),
\quad
\bfn_\gamma'(t) = - \hat v_\gamma(t) \bft_\gamma(t) \]
where $\hat v_\gamma(t) = \kappa_\gamma(t) v_\gamma(t)$,
$\kappa_\gamma(t)$ being the geodesic curvature of $\gamma$.
The curve $\gamma$ is locally convex
if and only if $\hat v_\gamma(t) > 0$ for all $t$.
In matrix notation,
$\gamma(t)$, $\bft_\gamma(t)$ and $\bfn_\gamma(t)$
are the columns of the orthogonal matrix $\fF_\gamma(t)$
which satisfies
\[ \fF_\gamma'(t) = \fF_\gamma(t) \Lambda_\gamma(t); \quad 
\Lambda_\gamma(t) = \begin{pmatrix} 0 & -v_\gamma(t) & 0 \\
v_\gamma(t) & 0 & -\hat v_\gamma(t) \\
0 & \hat v_\gamma(t) & 0 \end{pmatrix}. \]
Let $V \subset so_3$ be the plane
(i.e., $2$-dimensional real vector space)
of matrices $M$ with $(M)_{31} = 0$;
$\Lambda_\gamma$ can be considered a function from $[0,1]$ to $V$.
Let $V_{\cI} \subset V$ be the half-plane $(M)_{21} > 0$;
$\Lambda_\gamma$ can also be considered a function
$\Lambda_\gamma: [0,1] \to V_{\cI}$.
Conversely, given a smooth function $\Lambda_\gamma: [0,1] \to V_{\cI}$
or, equivalently, $v_\gamma$ and $\hat v_\gamma$,
the above equations together with $\fF_\gamma(0) = I$
may be interpreted as an initial value problem
defining $\fF_\gamma(t)$ and therefore $\gamma(t) = \fF_\gamma(t) e_1$.

Our aim is to consider a reasonably large space of functions
which still allows the initial value problem above to be solved.
The Hilbert space $L^2([0,1])$ is now a natural choice:
for $v, \hat v \in L^2([0,1])$ the initial value problem can be solved.
We therefore interpret $\cI_Q$ to be the closed subset of $(L^2([0,1]))^2$
of pairs of functions $(w, \hat v)$ such that
the solution $\Gamma: [0,1] \to SO_3$
of the initial value problem
\begin{equation} \label{eq:Gamma}
\Gamma'(t) = \Gamma(t) \Lambda(t), \quad
\Gamma(0) = I,
\end{equation}
satisfies $\Gamma(1) = Q$,
where
\[
\Lambda(t) = \begin{pmatrix} 0 & -v(t) & 0 \\
v(t) & 0 & -\hat v(t) \\
0 & \hat v(t) & 0 \end{pmatrix},
\quad
v(t) = \frac{w(t) + \sqrt{(w(t))^2 + 4}}{2}
\]
(we have $w(t) = v(t) - 1/v(t)$, $v(t) > 0$).
The fact that $\Lambda(t)$ belongs to the Lie algebra $so_3$
guarantees that $\Gamma(t)$ assumes values
in the corresponding Lie group $SO_3$.
Also, the function $\Gamma$ defined above is continuous,
absolutely continuous and differentiable almost everywhere.
In this way $\cI_Q \subset (L^2([0,1]))^2$
is a smooth Hilbert manifold of codimension $3$.
Indeed, consider the map
$\omega_\cI: (L^2([0,1]))^2 \to SO_3$ taking $(w, \hat v)$ to $\Gamma(1)$,
where $\Gamma: [0,1] \to SO_3$ is defined
by the initial value problem (\ref{eq:Gamma}) above.
Smooth dependence on parameters tells us that this map
is smooth; let us compute its derivative.

In general, for a curve $\Gamma: [0,1] \to SO_3$,
write $\Gamma(t_0;t_1) = (\Gamma(t_0))^{-1} \Gamma(t_1)$.
Let $L$ be a one-parameter family of functions
$L(s): [0,1] \to V_{\cI}$, $s \in (-\epsilon, \epsilon)$
with $L(0) = \Lambda$.
Let $G(s): [0,1] \to SO_3$ be the solution of
\[ (G(s))'(t) = (G(s))(t) (L(s))(t), \quad
(G(s))(0) = I \]
so that $G(0) = \Gamma$;
notice that the derivative in this initial value problem 
is with respect to $t$.
An easy computation gives that the derivative of $G$
(with respect to $s$) satisfies
\[ (G'(s))(t) (G(s)(t))^{-1} =
\int_0^t (G(s))(\tau) (L'(s))(\tau) ((G(s))(\tau))^{-1} d\tau. \]
Assume, for instance, that $(L'(0))(t)$ consists of three smooth
narrow bumps around times $t_i$ so that
\[ (\Gamma(1))^{-1} (G'(0))(t) \approx
\sum_i (\Gamma(t_i;1))^{-1} (L'(0))(t_i) \Gamma(t_i;1). \]
An easy computation shows that the spaces
$(\Gamma(t_i;1))^{-1} V \Gamma(t_i;1) \subset so_3$
are not constant and therefore $(\Gamma(1))^{-1} (G'(0))(t)$
may assume any value in $so_3$.
The derivative of $\omega_\cI$ is therefore surjective.
The map $\omega_\cI$ is thus a submersion and $\cI_Q$ a regular level set.
The geometric description of $\cI_Q$ comes from the identification
$(w,\hat v) \leftrightarrow \gamma$, where $\gamma(t) = \Gamma(t) e_1$.

Similarly, let $V_{\cL} \subset V \subset so_3$ be the quarter-plane
$(M)_{31} = 0$, $(M)_{32} > 0$, $(M)_{21} > 0$.
If $\gamma$ is locally convex then the image of $\Lambda_\gamma$
is contained in $V_{\cL}$ and, conversely,
given a smooth function $\Lambda: [0,1] \to V_{\cL}$,
equation (\ref{eq:Gamma}) obtains a locally convex curve.
Define $\cL_Q \subset (L^2([0,1]))^2$
to be the set of pairs $(w, \hat w)$ such that
$(w, \hat v) \in \cI_Q$ where
\[ \hat v(t) = \frac{\hat w(t) + \sqrt{(\hat w(t))^2 + 4}}{2}. \]
As above, define $\omega: (L^2([0,1]))^2 \to SO_3$
taking $(w, \hat w)$ to $\Gamma(1)$,
where $\Gamma$ is again defined by equation (\ref{eq:Gamma}).
The computation above also shows that the smooth map $\omega$
is a submersion 
and $\cL_Q$ is a smooth Hilbert submanifold of codimension $3$
in $(L^2([0,1]))^2$.
We have a natural injective map
$\cL_Q \to \cI_Q$ taking $(w,\hat w)$ to $(w, \hat v)$.
In the spirit of considering $\cL_Q$ and $\cI_Q$
to be sets of curves we call this map an inclusion
(even though it is not an isometry with the above metric).

The space $\cL_Q$ we just defined is rather large.
The space $\cL_Q^{[C^k]}$ of locally convex curves
of class $C^k$ ($k \ge 2$) is now naturally identified
to a dense subset $\cL_Q^{[C^k]} \subset \cL_Q$
with a different topology.
But are these two spaces similar?
Or, more precisely, is the inclusion a weak homotopy equivalence?
The answer here is yes.

In order to see this, first consider $\cL_Q^{[[C^k]]}$, $k \ge 0$,
the subset of $(C^k([0,1]))^2$ of pairs $(w, \hat w)$ such that
$\Gamma(1) = Q$
(where $\Gamma$ is defined by the initial value problem (\ref{eq:Gamma})).
The inclusion $\cL_Q^{[[C^k]]} \subset \cL_Q$ 
is a homotopy equivalence:
this follows directly from Theorem 2 from \cite{BST}.
For the convenience of the reader, we quote here
a simplified version of that result.

\begin{prop}
\label{prop:bstth2}
Let $X$ and $Y$ be separable Banach spaces.
Suppose $i:Y \to X$ is a bounded, injective linear map with dense image and 
$M \subset X$ a smooth, closed submanifold of finite codimension.
Then $N = i^{-1}(M)$ is a smooth closed submanifold of $Y$ and 
the restriction $i: N \to M$ is a homotopy equivalence.
\end{prop}

For $k \ge 1$ the curves $\gamma$ are now of class $C^2$
and we may assume them to be parametrized by a constant
multiple of arc length
(this does not change the homotopy type of the space
since the group of orientation-preserving
diffeomorphisms of $[0,1]$ is contractible).
In terms of the pairs $(w,\hat w)$,
this says that we may assume $w$ to be constant.
But for $w$ constant, $\hat w$ is of class $C^k$
if and only if $\gamma$ is of class $C^{k+2}$,
completing the argument.

Summing up, and simplifying this discussion a little,
we may assume our curves $\gamma$
to be as smooth as our constructions require
but when constructing a homotopy we may 
use curves for which the curvature is only piecewise continuous.
We will however try to be as consistent as reasonably possible
in the use of the large space $\cL$ defined above.

\section{Projective transformations and Bruhat cells}

\label{sect:bruhat}

In this section we present some more algebraic notation,
especially the decomposition of $SO_3$ and $\Ss^3$
in (signed) \textit{Bruhat cells}.
Some of this material is presented in \cite{SaSha}
in a more algebraic fashion and for arbitrary dimension.

The \textit{projective transformation}
$\pi(A): \Ss^2 \to \Ss^2$ associated
to $A \in SL_3$ is defined by
$\pi(A)(v) = \widehat{Av} = (1/|Av|) Av$.
Similarly, define $\pi(A): SO_3 \to SO_3$
by $\pi(A)(Q_0) = Q_1$ if there exists
$U_1 \in \Up$ with $AQ_0 = Q_1U_1$;
here $\Up$ is the contractible group of
upper triangular matrices with positive diagonal
and determinant $+1$.
In other words, $\pi(A)(Q_0)$ is obtained from $AQ_0$
by performing Gram-Schmidt on columns.
In particular, if $A \in SO_3$ then $\pi(A)(Q_0) = AQ_0$.
This action of $SL_3$ on $SO_3$ is transitive
but not doubly transitive;
we shall soon discuss the extent to which it fails
to be doubly transitive.

Notice that $\pi(A): \Ss^2 \to \Ss^2$ is
an orientation preserving diffeomorphism.
For an immersion $\gamma: [0,1] \to \Ss^2$
define $\pi(A) \gamma$ by $(\pi(A)\gamma)(t) = \pi(A)(\gamma(t))$;
the curve $\pi(A)\gamma$ is again an immersion.
Since $\pi(A)$ takes geodesics to geodesics,
if $\gamma: [0,1] \to \Ss^2$ is locally convex then so is $\pi(A) \gamma$
(this can also be checked directly from the definition
by a straightforward computation).
The map $\pi(A)$ is defined so that
$\fF_{\pi(A)\gamma} = \pi(A)(\fF_{\gamma})$
for any immersion $\gamma: [0,1] \to \Ss^2$.
Notice that $\pi(A)I = I$ if and only if $A \in \Up$.
Thus, for $U \in \Up$, $\pi(U)$ defines a smooth
map from $\cL_Q$ to $\cL_{\pi(U) Q}$.

Let $B_{3} \subset O_{3}$ be the Coxeter-Weyl group
of signed permutation matrices:
the group $B_3$ has $48$ elements
and corresponds to the isometries of the octahedron of vertices
$\pm e_1, \pm e_2, \pm e_3$.
Let $B^{+}_{3} = B_{3} \cap SO_{3}$ be the subgroup
of orientation preserving isometries.
Dropping signs defines a homomorphism
from $B^{+}_3$ to the symmetric group $S_3$:
the \textit{number of inversions} of $Q \in B^{+}_3$
is the number of inversions of the corresponding permutation in $S_3$
(recall that the number of inversions of $\pi \in S_n$
is the number of pairs $(i,j)$ with $i < j$ and $\pi(i) > \pi(j)$).
We denote an element of $B^{+}_3$ by the letter $P$
(for permutation), indicating in the subscript
the corresponding permutation (as a cycle) and the signs, read by column,
in binary:
\[ 0 = +++,\; 1 = ++-, \;2 = +-+, \;3 = +--, \;\ldots, \;7 = ---; \]
thus, for instance
\[
P_{(13);1} = \begin{pmatrix}
0 & 0 & -1 \\ 0 & +1 & 0 \\ +1 & 0 & 0 
\end{pmatrix}, \quad
P_{(13);2} = \begin{pmatrix}
0 & 0 & +1 \\ 0 & -1 & 0 \\ +1 & 0 & 0 
\end{pmatrix}. \]
The lift $\tilde B^{+}_3 \subset \Ss^3$ of $B^{+}_3$
is the group of $48$ quaternions with either
one coordinate of absolute value $1$ and three equal to $0$ or
two coordinates of absolute value $1/\sqrt{2}$ and two equal to $0$ or
four coordinates of absolute value $1/2$.
As a subset of $\RR^4$, $\tilde B^{+}_3$ is also the root system $F_4$
(with all vectors of size $1$, unlike what is usual for Lie algebras).
For instance,
\[ \Pi\left(\frac{\bfone-\bfj}{\sqrt{2}}\right) = P_{(13);1}, \quad
\Pi\left(\frac{\bfi+\bfk}{\sqrt{2}}\right) = P_{(13);2}. \]

The Bruhat cell $\Bruhat(Q) \subset SO_3$ of $Q \in SO_3$
is the set of all orthogonal matrices of the form
$U_0 Q U_1^{-1}$, with $U_0, U_1 \in \Up$.
Each Bruhat cell contains a unique element of $B^{+}_3$.
We also denote the Bruhat cells by $\Bruhat_{\ast}$,
with subscripts defined as for $P_{\ast} \in B^{+}_3$;
thus, for instance, $\Bruhat(P_{(13);1}) = \Bruhat_{(13);1}$.
The cell $\Bruhat(P)$ ($P \in B^{+}_3$) is diffeomorphic to $\RR^d$,
where $d$ is the number of inversions of $P$.
There are therefore exactly four open cells,
corresponding to the two matrices above plus
\[ \Pi\left(\frac{\bfone+\bfj}{\sqrt{2}}\right) = P_{(13);4}, \quad
\Pi\left(\frac{\bfi-\bfk}{\sqrt{2}}\right) = P_{(13);7}. \]

If $Q_0$ and $Q_1$ belong to the same Bruhat cell
then there exists a canonical and explicit diffeomorphism between
the Hilbert manifolds $\cL_{Q_0}$ and $\cL_{Q_1}$.
Indeed, for $Q_0$ and $Q_1 \in \Bruhat(Q_0)$
let $U_0, U_1 \in \Up$ be such that $U_0 Q_0 = Q_1 U_1$;
$U_0$ is uniquely defined if we require $U_0 \in \Upum$,
where $\Upum \subset \Up$ is the subgroup of matrices
with diagonal entries equal to $+1$.
Then the map $\pi(U_0):  \cL_{Q_0} \to \cL_{Q_1}$ is a diffeomorphism;
its inverse is the similarly constructed map $\pi(U_0^{-1})$.

\begin{figure}[ht]
\begin{center}
\epsfig{height=30mm,file=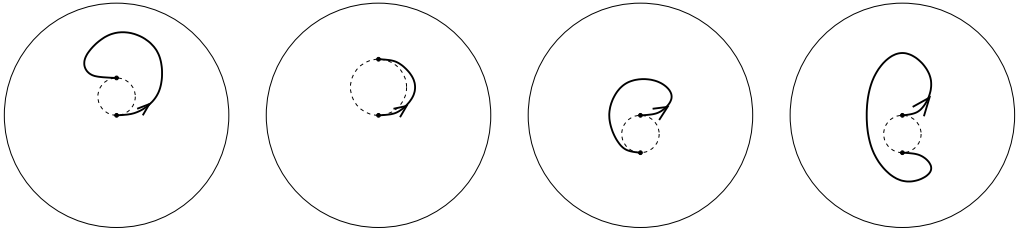}
\end{center}
\caption{Representatives of the open Bruhat cells}
\label{fig:genericbruhat}
\end{figure}

Figure \ref{fig:genericbruhat} shows example of curves
$\gamma \in \cL_Q$ for $Q \in \Bruhat_{(13);\ell}$
for $\ell = 1, 2, 4, 7$ (in this order).

The dashed closed convex curves indicate a convenient
geometric way to recognize these Bruhat cells:
in all cases there exists a closed convex curve
tangent to both endpoints of $\gamma$
and orientations at the endpoints allow
us to distinguish between cells.

For $\gamma \in \cL$,
let $\Gamma = \fF_\gamma: [0,1] \to SO_3$
and write $\fF_\gamma(t_0; t) = \Gamma(t_0; t) =
(\Gamma(t_0))^{-1} \Gamma(t)$.
Similarly, $\tilde\fF_\gamma(t_0;t) = \tilde\Gamma(t_0; t) =
(\tilde\Gamma(t_0))^{-1} \tilde\Gamma(t)$
where $\tilde\Gamma: [0,1] \to \Ss^3$ is a lift of $\Gamma$.
Clearly, $\Gamma(0; t) = \Gamma(t)$,
$\Gamma(t_0; t_0) = I$, $\tilde\Gamma(t_0; t_0) = \bfone$.
On the other hand, there exists $\epsilon > 0$
such that for all $t \in (t_0, t_0 + \epsilon)$
we have $\Gamma(t_0; t) \in \Bruhat_{(13),2}$.

\begin{figure}[ht]
\begin{center}
\epsfig{height=20mm,file=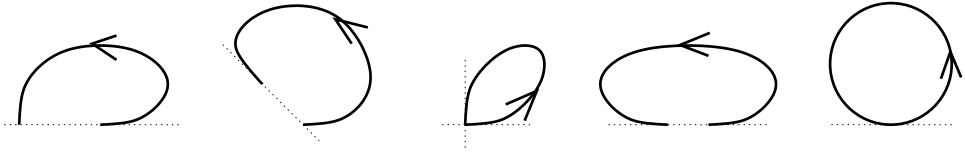}
\end{center}
\caption{Convex arcs}
\label{fig:5bruhat}
\end{figure}

A locally convex arc $\gamma|_{[t_0,t_1]}$
is \textit{convex} if
$\fF_\gamma(t_0;t) \in \Bruhat_{(13);2}$
for all $t \in (t_0,t_1)$.
Notice that we do not require that
$\fF_\gamma(t_0;t_1)  \in \Bruhat_{(13);2}$;
if this also happens then $\gamma$ is \textit{stably convex}.
There are five other Bruhat cells
to which $\fF_\gamma(t_0;t_1)$ may belong:
$\Bruhat_{(123);6}$, $\Bruhat_{(132);0}$,
$\Bruhat_{(23);2}$, $\Bruhat_{(12);4}$
and $\Bruhat_{e;0} = I$.
Examples of convex arcs corresponding
to these five cells are given in Figure \ref{fig:5bruhat}.
We shall come back to these five cells again and again.

Recall that a matrix $Q \in SO_3$ is (stably) \textit{convex}
if there exists a (stably) convex arc $\gamma: [t_0, t_1] \to \Ss^2$
with $\fF_\gamma(t_0;t_1) = Q$.
Thus, the set of stably convex matrices is the open cell 
$\Bruhat_{(13);2}$ and the set of convex matrices is the union
of the six Bruhat cells above (including the open cell).
Also, a quaternion $z \in \Ss^3$ is (stably) \textit{convex}
if there exists a (stably) convex arc $\gamma: [t_0, t_1] \to \Ss^2$
with $\tilde\fF_\gamma(t_0;t_1) = z$.
The six convex elements of $\tilde B_3^{+}$ are
\[ \frac{\bfi+\bfk}{\sqrt{2}}, \quad 
\frac{-\bfone+\bfi-\bfj+\bfk}{2}, \quad
\frac{-\bfone+\bfi+\bfj+\bfk}{2}, \quad
\frac{-\bfone+\bfi}{\sqrt{2}}, \quad
\frac{-\bfone+\bfk}{\sqrt{2}}, \quad
-\bfone; \]
only the first one is stably convex.
A quaternion $z$ is (stably) \textit{anticonvex}
if $-z$ is (stably) convex.
The sets of convex and anticonvex quaternions
are unions of distinct Bruhat cells and therefore disjoint;
furthermore, the only points of intersection
between their respective closures are $\pm\bfone$.

\medskip

Let $\cC_{\nu}$ be the circle (contained in $\Ss^2$) with diameter $e_1e_3$,
parametrized by $\nu_1 \in \cL_I$,
\[ \nu_1(t) = \frac{1}{2}
\left( {1 + \cos(2 \pi t)},
{\sqrt{2}}\; \sin(2 \pi t),
{1 - \cos(2 \pi t)} \right) \]
for which
\[ \tilde\fF_{\nu_1}(t) =
\exp\left(\pi t\,\bfkh\right), \quad
\tilde\Lambda_{\nu_1}(t) = \pi\bfkh, \quad
\bfkh = \frac{\bfi+\bfk}{\sqrt 2},
\]
where $\tilde\Lambda_{\gamma}(t) =
(\tilde\fF_\gamma(t))^{-1}\tilde\fF'_\gamma(t)$.
For a positive real number $s$, let $\nu_s(t) = \nu_1(st)$
so that $\nu_s \in \cL_{\exp(\pi s \bfkh)}$.
In particular,
for integer $k > 0$, $\nu_k \in \cL_{(-\bfone)^k}$.
We also have $\nu_1 \in \cL_{-\bfone,c}$ and
$\nu_k \in \cL_{-\bfone,n}$ for $k$ odd, $k > 1$.

\medskip

More generally, a circle of radius $\rho < \pi/2$ is a closed convex curve:
\[ \gamma(t) =
\cos(2\pi t) \sin(\rho) v_1 + 
\sin(2\pi t) \sin(\rho) v_2 + 
\cos(\rho) v_3, \quad t \in [0,1], \]
where $v_1, v_2, v_3$ is a positively oriented orthonormal basis.
Thus, $\nu_s$ is a circle of radius $\pi/4$
(measured along the sphere).

\medskip

The image of a convex circle by a projective transformation
is a \textit{spherical ellipse}, or just \textit{ellipse}.
Notice that for us an ellipse is an {oriented} curve.
Also, a \textit{projective arc-length parametrization} of an ellipse
is a locally convex curve $\gamma: [t_0, t_1] \to \Ss^2$
of the form $\gamma = \pi(A) \circ \tilde\gamma$
where $A \in SL_3$ and $\tilde\gamma$ is a parametrization
by a multiple of arc length of a circle.
We shall sometimes use ellipses
when a concrete choice of convex arc is desirable.

\begin{lemma}
\label{lemma:ellipse}
Let $\gamma: [0,1] \to \Ss^2$ be a stably convex arc
and let $t \in (0,1)$.
Set 
\[ \tilde\fF_\gamma(0) = z_0, \quad \gamma(t) = v_{t},
\quad \tilde\fF_\gamma(1) = z_1. \]
Then there exists $\epsilon > 0$ such that
if $\hat z_0 \in \Ss^3$, $\hat v_{t} \in \Ss^2$ and $\hat z_1 \in \Ss^3$
are such that
\[ |\hat z_0 - z_0| < \epsilon, \quad
|\hat v_{t} - v_{t}| < \epsilon, \quad
|\hat z_1 - z_1| < \epsilon \]
then there exists a unique ellipse $\hat\cE \subset \Ss^2$
and projective arc-length parametrization
$\hat\gamma: [0,1] \to \hat\cE \subset \Ss^2$ with
\[ \tilde\fF_{\hat\gamma}(0) = \hat z_0, \quad \hat\gamma(t) = \hat v_{t},
\quad \tilde\fF_{\hat\gamma}(1) = \hat z_1. \]
\end{lemma}

{\nobf Proof:}
Let $Q_i = \Pi(z_i)$;
draw in $\Ss^2$ geodesics $\ell_i$
perpendicular to $Q_ie_3$, oriented by $Q_ie_2$ at $Q_ie_1$
so that $\ell_i$ is tangent to $\gamma$ at $t = i$ ($i=0,1$).
Since $z_0^{-1}z_1 \in \Bruhat_{(13);2}$,
the geodesics $\ell_0$ and $\ell_1$ are transversal
and divide the sphere into four open regions,
with the image of $\gamma$ contained in the region
characterized by being to the left of both $\ell_0$ and $\ell_1$.
Since all the relevant conditions are open,
for sufficiently small $\epsilon > 0$
the corresponding geodesics $\hat\ell_i$
are also transversal and $\hat v_{t}$ lies to the left
of both $\hat\ell_0$ and $\hat\ell_1$.
There exists $A \in SL_3$ such that
the projective transformation $\pi(A)$ satisfies
$\pi(A)\hat Q_0 = I$ and $\pi(A)\hat Q_1 = P_{(13);2}$;
for sufficiently small $\epsilon > 0$
we may furthermore assume that $\pi(A) \hat v_{t}$
lies in the first octant.
There is $c > 0$ such that, for
\[ B_c = \begin{pmatrix} c & 0 & 0 \\
0 & c^{-2} & 0 \\ 0 & 0 & c \end{pmatrix}, \]
$\pi(B_cA) \hat v_{t}$ lies on the arc of circle $\nu_1$ defined above,
thus obtaining the desired arc of ellipse.
Uniqueness follows from the fact that five points determine a conic
and three points determine a projective transformation in the line.
\qed

Recall that two smooth curves \textit{osculate} each
other at a common point if they are tangent and have the same curvature at that point.
The next lemma may be considered a limit case of the previous one.

\begin{lemma}
\label{lemma:osculate}
Let $Q \in \Bruhat_{(13);2}$;
then there exists a unique ellipse $\cE \subset \Ss^2$
and projective arc-length parametrization
$\gamma: [0,1] \to \cE \subset \Ss^2$
with $\gamma$ convex,
$\gamma(0) = e_1$, $\gamma'(0) = e_2$,
$\cE$ osculating the circle $\cC_{\bfkh}$ at $e_1$
and $\fF_\gamma(1) = Q$.
\end{lemma}

{\nobf Proof:} As in the previous proof,
there exists $A \in SL_3$ with $\pi(A)I = I$
and $\pi(A)Q = P_{(13);2}$.
With $B_c$ as in the previous proof,
the ellipses $\gamma_c = \pi(A^{-1}B_c) \circ \nu_1$
satisfy $\gamma_c(0) = e_1$, $\gamma_c'(0) = e_2$,
and $\fF_{\gamma_c}(1) = Q$;
there exists a unique $c > 0$ for which $\gamma_c$ oscullates $\cC_{\bfkh}$.
Uniqueness again follows from five points determining a conic
(oscullation counts as three points).
\qed


\section{Total curvature}

\label{sect:total}

Given a locally convex curve $\gamma: [a,b] \to \Ss^2$,
let $\kappa_\gamma: [a,b] \to \RR$ be the geodesic curvature of $\gamma$.
This must not be confused with the euclidean curvature
$\kappa^E_\gamma$ of $\gamma$ interpreted
as a curve in $\RR^3$:
$\kappa^E_\gamma(t) = \sqrt{1 + \kappa_\gamma^2(t)}$.
The total (euclidean) curvature of $\gamma: [a,b] \to \Ss^2$ is
\[ \tot(\gamma) = \int_{[a,b]} \kappa^E_\gamma(t)\; |\gamma'(t)| dt. \]
Notice that $\tot(\nu_s) = 2\pi s$.

\begin{lemma}
\label{lemma:total}
Let $\gamma: [a,b] \to \Ss^2$ be a locally convex curve.
\begin{enumerate}[(a)]
\item{The total curvature of $\gamma$ equals
the total variation of $\bft_\gamma$ and
twice the length of $\tilde\fF_\gamma$:
\[ \tot(\gamma) = \int_{[a,b]} |\bft'_\gamma(t)|\;dt
= 2 \int_{[a,b]} |\tilde\fF_\gamma'(t)|\;dt. \]}
\item{If $\gamma \in \cL_I$ is a closed convex curve
then $\tot(\gamma) \in [2\pi, 4\pi)$.}
\item{If $\gamma \in \cL$ is a convex arc
then $\tot(\gamma) \in (0, 4\pi]$.}
\item{If $\tot(\gamma) < \pi$ then $\gamma$ is convex.}
\end{enumerate}
\end{lemma}

The inequalities above are not necessarily the best possible.

{\nobf Proof:}
Item (a) is a straightforward computation.
The first inequality in item (b)
is the well known general fact that the total
curvature of a closed curve is at least $2\pi$.
Alternatively, given (a),
$\tilde\Gamma = \tilde\fF_\gamma: [0,1] \to \Ss^3$
satisfies $\tilde\Gamma(0) = \bfone$, $\tilde\Gamma(1) = -\bfone$
so of course the length of $\tilde\Gamma$ must be at least $\pi$.
For the second inequality in (b), first notice that
$\kappa^E_\gamma(t) \le 1 + \kappa_\gamma(t)$ so we must prove that
\[ I_1 + I_2 < 4\pi, \quad
I_1 = \int_{[0,1]} |\gamma'(t)| dt,
\quad
I_2 = \int_{[0,1]} \kappa_\gamma(t)\;  |\gamma'(t)| dt. \]
By Gauss-Bonnet, $I_2 = 2\pi - A_\gamma < 2\pi$,
where $A_\gamma$ is the area of the smaller disk in $\Ss^2$
bounded by $\gamma$.
Let $\bfn_\gamma: {[0,1]} \to \Ss^2$ be the unit normal vector:
$\bfn_\gamma(t) = \fF_\gamma(t) e_3$.
A straightforward computation shows that
\[ |\bfn_\gamma'(t)| = \kappa_\gamma(t) |\gamma'(t)|, \quad
\kappa_{\bfn_\gamma}(t) = 1/\kappa_\gamma(t). \]
Thus, again by Gauss-Bonnet,
\[ I_1 = \int_{[0,1]} \kappa_{\bfn_\gamma}(t) \;  |\bfn_\gamma'(t)| dt =
2\pi - A_{\bfn_\gamma} < 2\pi \]
where $A_{\bfn_\gamma}$ is (of course) the area of the smaller
disk in $\Ss^2$ bounded by $\bfn_\gamma$.
This completes the proof of (b).

For item (c), consider a convex arc $\gamma:[0,1] \to \Ss^2$.
For arbitrarily small $\epsilon > 0$,
$\gamma|_{[0,1-\epsilon]}$ can be extended
to a closed convex curve $\tilde\gamma$.
Thus, from (b), $\tot(\gamma|_{[0,1-\epsilon]}) < \tot(\tilde\gamma) < 4\pi$.
Since this estimate holds for all $\epsilon$,
$\tot(\gamma) \le 4\pi$.

For item (d), assume that $\gamma_1: [a,b] \to \Ss^2$ is convex but that
if the domain of $\gamma_1$ is extended to define $\gamma_2$
then $\gamma_2$ is not convex.
In other words, $\gamma_1$ is similar to one of the curves
in Figure \ref{fig:5bruhat}.
A case by case analysis shows that there exist $t_0 < t_1 \in [a,b]$
with $\bft_{\gamma_1}(t_1) = -\bft_{\gamma_1}(t_0)$
and therefore 
\[ \tot(\gamma_1) \ge \int_{[t_0,t_1]} |\bft_{\gamma_1}'|(t)\;dt \ge \pi. \]
\qed

\section{The map $g_0$}

\label{sect:2to4}

The aim of this section is to construct a map
$g_0: \Ss^2 \to \cL_{\bfone} \subset \cI_{\bfone}$.
We shall have $g_0(\bfs) = \nu_2$ and $g_0(\bfn) = \nu_4$,
where $\bfn = (0,0,+1)$ and $\bfs = (0,0,-1)$
are the north and south pole, respectively.
The existence of such a map proves
that $\nu_2$ and $\nu_4$ are in the same connected component of $\cL_{\bfone}$,
consistently with Little's result that $\cL_I$ has three connected components 
$\cL_{-\bfone, c}$, $\cL_{\bfone}$ and $\cL_{-\bfone,n}$
(see Figure \ref{fig:3comp}).
As we shall see in Proposition \ref{prop:g0generator},
$g_0$ turns out to be a generator
of $\pi_2(\cI_{\bfone}) \approx H_2(\cI_{\bfone};\ZZ) \approx \ZZ$.
The map $g_0$ is one the the crucial objects in this paper
so its construction shall be discussed in some detail.

\begin{figure}[hp]
\begin{center}
\epsfig{height=190mm,file=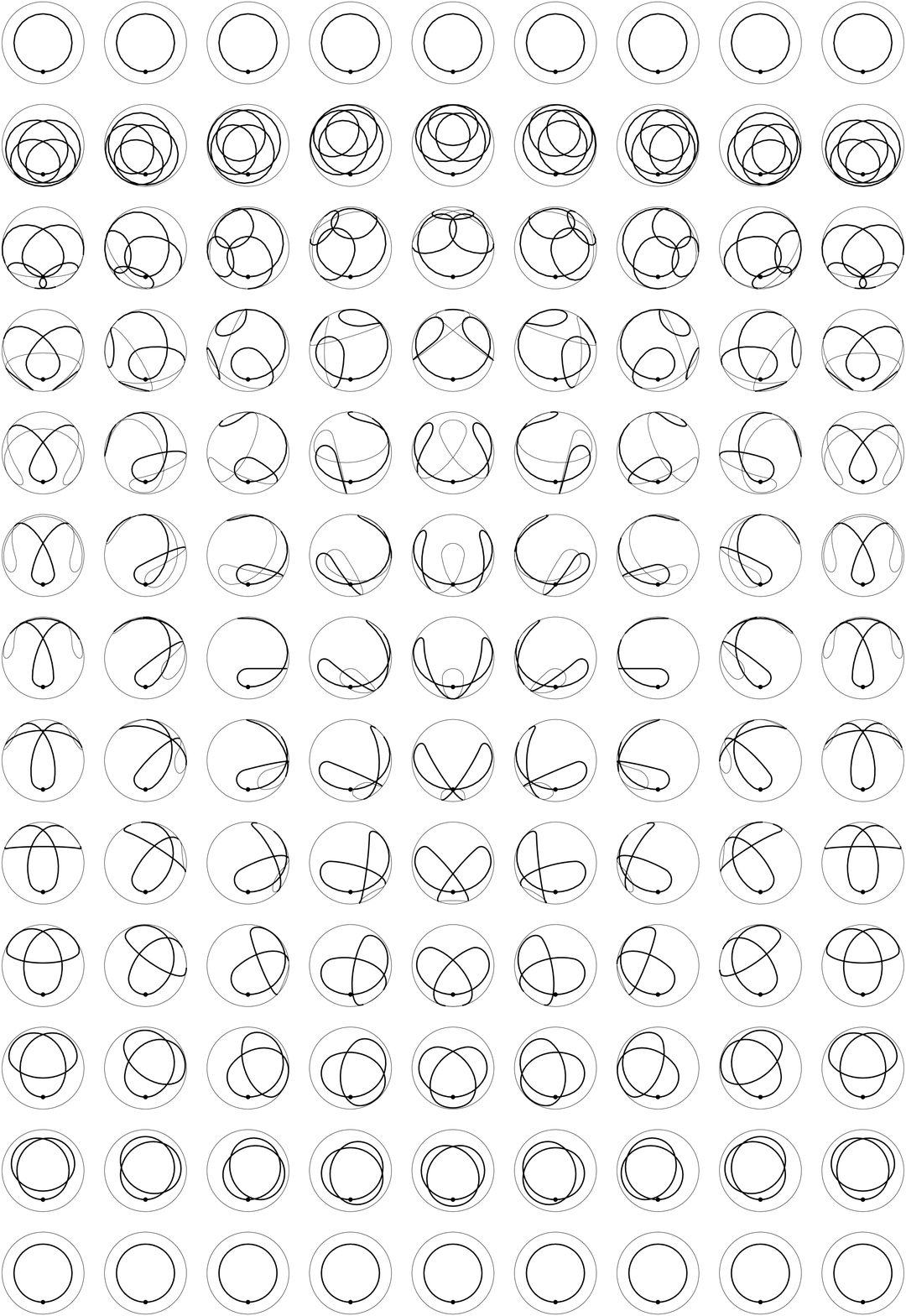}
\end{center}
\caption{A sketch of $g_0: \Ss^2 \to \cL_{\bfone}$.}
\label{fig:Gamma5}
\end{figure}

The map $g_0$ is shown in Figure \ref{fig:Gamma5}:
here the bottom line is the south pole $\bfs = (0,0,-1)$,
the top line is the north pole $\bfn = (0,0,+1)$
and other horizontal lines are circles
contained in a plane of the form $z = z_0$;
vertical lines are meridians, i.e., half circles
joining the two poles.
Each small sphere is drawn here as a photo of a transparent sphere:
a curve in the front is drawn a bit thicker.
The vector perpendicular to the page pointing towards the reader
is $e_1+e_3$ with $e_3 - e_1$ pointing up and $e_2$ to the right of the reader.
The base point $e_1$ is drawn as a thick dot.

Consider the equator of the unit sphere $\Ss^2$
contained in the plane $z = 0$.
Fix a unit normal vector $N = (0,0,1)$ to the plane and call it \textit{up}.
Consider six equally spaced points
$P_0 = P_3 = (1/2,-\sqrt{3}/2,0)$, $Q_0 = Q_3 = (1,0,0)$,
$P_1 = (1/2,\sqrt{3}/2,0)$, $Q_1 = -P_0$, $P_2 = -Q_0$ and $Q_2 = -P_1$
along the equator.
Notice that $Q_i/2$ is the midpoint of the segment $P_iP_{i+1}$.
For $\alpha \in \RR$, let $\tilde\alpha = \arcsin(\sin(\alpha)/2)$
so that $\tilde\alpha \in [-\pi/6,\pi/6]$.
Let 
\[ Q_i^{\pm}(\alpha) =
\cos(\alpha-\tilde\alpha) Q_i \pm \sin(\alpha-\tilde\alpha) N, \]
so that $Q_i^{+}(0) = Q_i^{-}(0) = Q_i$, $Q_i^{-}(\alpha) = Q_i^{+}(-\alpha)$
and $Q_i^{+}$ (a function of $\alpha$)
parametrizes the circle passing through $\pm Q_i$ and $\pm N$.
Equivalently, $Q_i^{+}(\alpha)$ is the only point on the unit sphere
such that the vector $Q_i^{+}(\alpha) - (Q_i/2)$
is a positive multiple of $(\cos\alpha)Q_i + (\sin\alpha)N$.
Let $\cC_i^{\pm}(\alpha)$ be the circle containing
$P_i$, $Q_i^{\pm}(\alpha)$ and $P_{i+1}$,
so that $\cC_i^{\pm}(\alpha)$ has radius $\cos(\tilde\alpha)$
and is contained in a plane perpendicular to 
$- \cos(\alpha) N \pm \sin(\alpha) Q_i$.
Let $A_i^{\pm}(\alpha) \subset \cC_i^{\pm}(\alpha)$
be the arc from $P_i$ through $Q_i^{\pm}(\alpha)$ to $P_{i+1}$.

\begin{figure}[ht]
\psfrag{x}{$x$}
\psfrag{y}{$y$}
\psfrag{N}{$N$}
\psfrag{P0}{$P_0$}
\psfrag{P1}{$P_1$}
\psfrag{P2}{$P_2$}
\psfrag{Q0}{$Q_0$}
\psfrag{Q1}{$Q_1$}
\psfrag{Q2}{$Q_2$}
\psfrag{aaf}{$\alpha$}
\psfrag{atil}{$\tilde\alpha$}
\psfrag{Qi}{$Q_i$}
\psfrag{Qiplus}{$Q_i^{+}(\alpha)$}
\psfrag{PiPiplus}{$\frac{Q_i}{2}$}
\begin{center}
\epsfig{height=50mm,file=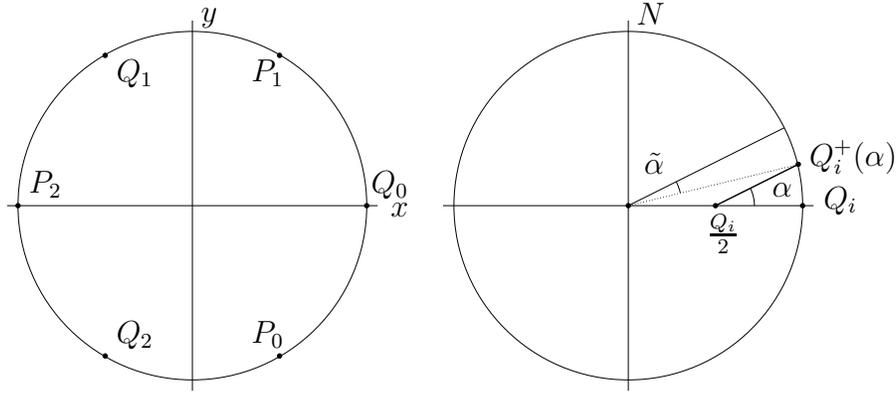}
\end{center}
\caption{Constructing $P_i$, $Q_i$ and $Q_i^{+}(\alpha)$.}
\label{fig:PiQi}
\end{figure}


Orient the circles $\cC_i^{\pm}(\alpha)$ and the arcs $A_i^{\pm}(\alpha)$
from $P_i$ through $Q_i^{\pm}(\alpha)$ to $P_{i+1}$.
The circle $\cC_i^{\pm}(\alpha)$ has geodesic curvature
equal to $\mp\tan(\tilde\alpha)$.
Parametrize the arcs 
$A_0^{+}(\alpha)$, $A_1^{-}(\alpha)$, $A_2^{+}(\alpha)$,
$A_0^{-}(\alpha)$, $A_1^{+}(\alpha)$, $A_2^{-}(\alpha)$ 
and $A_0^{+}(\alpha)$
by a multiple of arc length using the domains
$[-1/12,1/12]$, $[1/12,3/12]$, $[3/12,5/12]$,
$[5/12,7/12]$, $[7/12,9/12]$, $[9/12,11/12]$
and $[11/12,13/12]$, respectively.
Concatenate the above parametrizations to define
a parametrization $\beta_\alpha: [0,1] \to \Ss^2$
by a multiple of arc length of a curve $C_\alpha \subset \Ss^2$.
In particular, $\beta_\alpha(0) = Q_0^{+}(\alpha)$,
$\beta_\alpha(1/12) = P_1$, $\beta_\alpha(2/12) = Q_1^{-}(\alpha)$ and so on.
Notice that $\beta_\alpha$ is of class $C^1$,
even at the points $t = j/12$.
The curve $\beta_\alpha$ is an immersion;
its geodesic curvature is always in the interval $[-\sqrt{3}/3,+\sqrt{3}/3]$,
with the extremes assumed for $\alpha = \pi/2 + k\pi$, $k \in \ZZ$.
The curve $\beta_0$ is the equator covered twice;
$\beta_\pi$ is the equator covered four times with the opposite orientation.
Define $B_\alpha = \fF_{\beta_\alpha}: [0,1] \to SO_3$ as usual;
lift this to define $\tilde B_\alpha: [0,1] \to \Ss^3$
with $\tilde B_0(0) = \bfone$ and $\tilde B_\alpha(t)$
a continuous function of $\alpha$ and $t$,
$1$-periodic in $t$ and $4\pi$-periodic in $\alpha$.


Define $\bfh = \exp(\pi\bfj/8)$; notice that
\[ \bfh\bfi\bfh^{-1} = \bfih = \frac{\bfi-\bfk}{\sqrt{2}}, \quad
\bfh\bfk\bfh^{-1} = \bfkh = \frac{\bfi+\bfk}{\sqrt{2}}. \]
Define $\tilde\Gamma_\alpha(t) = \tilde B_\alpha(t) \bfh^{-1}$,
$\Gamma_\alpha = \Pi \circ \tilde\Gamma_\alpha$ and (of course)
\[ \gamma_\alpha(t) = \Gamma_\alpha(t) e_1 =
B_\alpha(t) \frac{e_1+e_3}{\sqrt{2}} =
\frac{\sqrt{2}}{2} \beta_\alpha(t) +
\frac{\sqrt{2}}{2} \bfn_{\beta_\alpha}(t). \]
Geometrically, $\gamma_\alpha$ is the concatenation
of six circle arcs $\tilde A^{\pm}_i$,
obtained from $A^{\pm}_i$ by increasing or decreasing
the radius by $\pi/4$.
The geodesic curvature of $\gamma_\alpha$ is
\[ \kappa_{\gamma_\alpha}(t) = \begin{cases}
\tan\left(\frac{\pi}{4}-\tilde\alpha\right), &
t \in [0,\frac{1}{12}) \cup (\frac{3}{12},\frac{5}{12}) \cup 
(\frac{7}{12},\frac{9}{12}) \cup (\frac{11}{12},1], \\
\tan\left(\frac{\pi}{4}+\tilde\alpha\right), &
t \in (\frac{1}{12},\frac{3}{12}) \cup 
(\frac{5}{12},\frac{7}{12}) \cup (\frac{9}{12},\frac{11}{12}).
\end{cases} \]
and therefore always in the interval $[2-\sqrt{3},2+\sqrt{3}]$,
with extreme values assumed for $\alpha = \pi/2$.
Furthermore, $\tot(\gamma_\alpha)$ is a strictly increasing 
function of $\alpha \in [0,\pi]$
with $\tot(\gamma_0) = 4\pi$, $\tot(\gamma_\pi) = 8\pi$.
Notice that
\[ \tilde\Gamma_\alpha\left(t+\frac{1}{3}\right) =
\exp\left(\frac{4\pi}{3}\bfk\right) \tilde\Gamma_\alpha(t);
\quad
\tilde\Gamma_{\alpha}\left(t+\frac{1}{2}\right) =
-\tilde\Gamma_{-\alpha}(t).  \]
Define $g_0: \Ss^2 \to \cL_{\bfone}$ by
\[
g_0(p)(t) =
\left( \Gamma_\alpha\left( \frac{\theta}{6\pi} \right) \right)^{-1}
\gamma_\alpha\left( t + \frac{\theta}{6\pi} \right),
\quad
p = (\cos\theta \sin\alpha, \sin\theta \sin\alpha, -\cos\alpha).
\]
Notice that $g_0$ is well defined.

This construction yields the following result,
due to Little (\cite{Little});
we state and prove it here in order to get used
to the above construction which shall be needed later.

\begin{lemma}
\label{lemma:2to4}
Let $n$ be a positive integer, $n > 1$.
The curves $\nu_n$ and $\nu_{n+2}$ are in the same connected component
of $\cL_{I}$.
\end{lemma}

{\nobf Proof:}
The case $n=2$ follows from the above construction of $g_0$.
For larger $n$, just consider $\nu_n$ as a concatenation
of $\nu_{n-2}$ with $\nu_2$,
keep $\nu_{n-2}$ fixed and apply the above construction
to move from $\nu_2$ to $\nu_4$,
thus obtaining a path in $\cL_I$ from $\nu_n$ to $\nu_{n+2}$.
\qed

We prove an auxiliary result concerning the curves
$\beta_\alpha$ and $\gamma_\alpha$ for later use.
A \textit{common oriented tangent} to two oriented circles
is an oriented geodesic which is tangent to both circles,
with compatible orientation at tangency points.
Thus, for instance, $\cC_i^{\pm}(\alpha)$ and $\cC_{i+1}^{\mp}(\alpha)$
are tangent at $P_{i+1}$;
they also have a common oriented tangent,
a geodesic passing through $P_{i+1}$.
Let $\cC_{\bfk}$ (resp. $\cC_{\bfkh}$) be the great circle
and subgroup of $\Ss^3$
of points of the form $\exp(s \bfk)$ (resp. $\exp(s \bfkh)$),
$s \in \RR$.
Recall that we write $\tilde B_\alpha(t_0; t_1) =
(\tilde B_\alpha(t_0))^{-1} \tilde B_\alpha(t_1)$.

\begin{lemma}
\label{lemma:nocommontangent}
Let $\alpha \in (0, \pi)$.
\begin{enumerate}[(a)]
\item{The common oriented tangents between two distinct
circles among $\cC_i^{\pm}(\alpha)$, $i = 0, 1, 2$,
are the geodesics passing through the tangency point $P_{i+1}$
between $\cC_i^{\pm}(\alpha)$ and $\cC_{i+1}^{\mp}(\alpha)$, and only these.}
\item{For $t_0, t_1 \in [0,1)$,
if $\tilde B_\alpha(t_0; t_1) \in \cC_{\bfk}$
then $t_0 = t_1$.}
\item{For $t_0, t_1 \in [0,1)$,
if $\tilde\Gamma_\alpha(t_0; t_1) \in \cC_{\bfkh}$
then $t_0 = t_1$.}
\end{enumerate}
\end{lemma}

{\nobf Proof:}
Recall that our circles $\cC_i^{\pm}(\alpha)$
are not geodesics.
Each circle therefore defines a disk 
(the smaller connected component of the complement).

For item (a),
we have three essentially different pairs of circles
to consider: $(\cC_0^{+}(\alpha),\cC_0^{-}(\alpha))$,
$(\cC_0^{+}(\alpha),\cC_1^{+}(\alpha))$ and
$(\cC_0^{+}(\alpha),\cC_1^{-}(\alpha))$.
In the first case, the two circles have opposite
orientations and therefore the circles and corresponding disks
would have to lie on opposite sides of a common oriented tangent.
Since the open disks intersect, no common oriented tangent exists.
In the second case orientations agree and therefore
both disks would lie on the same side of a common oriented tangent.
But the union of the two closed disks contains the arc
from $P_0$ through $P_1$ to $P_2$ and therefore
is not contained in a hemisphere.
Thus also in this case no common oriented tangent exists.
Finally, in the third case again the two circles have opposite
orientations and therefore the disks must lie on opposite sides
of a common oriented tangent.
But the closed disks touch at $P_1$:
the common oriented tangent must therefore pass through this point,
completing the proof of the first claim.

For item (b), assume by contradiction
that $\tilde B_\alpha(t_0; t_1) = \exp(s\bfk)$, $t_0 \ne t_1$.
Consider the curve $\tilde\Gamma: [0,1] \to \Ss^3$
given by $\tilde\Gamma(t) = \tilde B_\alpha(t_0) \exp(2\pi t \bfk)$.
Notice that $\tilde\Gamma(0) = \tilde B_\alpha(t_0)$
and $\tilde\Gamma(s/(2\pi)) = \tilde B_\alpha(t_1)$.
Also $\tilde\Gamma(t) = \tilde\fF_{\gamma}(t)$
for $\gamma$ the geodesic
$\gamma(t) = B_\alpha(t_0) (\cos(4\pi t), \sin(4\pi t), 0)$.
Thus $\gamma$ is an oriented tangent to the curve $C_\alpha$
at two distinct points.
These two points must belong to different circles $\cC_i^{\pm}$
for a circle can not be twice tangent to the same geodesic.
But this contradicts item (a).

Finally, for item (c), notice that
\[ \tilde \Gamma_\alpha(t_0; t_1) =
\bfh \tilde B_\alpha(t_0;t_1) \bfh^{-1} \]
and therefore $\tilde\Gamma_\alpha(t_0; t_1) = \exp(s\bfkh)$
implies 
\[ \tilde B_\alpha(t_0;t_1) = \bfh^{-1} \exp(s\bfkh) \bfh
= \exp\left( s \bfh^{-1} \bfkh \bfh \right) = \exp(s \bfk). \]
Thus item (b) implies item (c).
\qed

\begin{prop}
\label{prop:g0generator}
The map $g_0$ is a generator of $\pi_2(\cI_{\bfone})$.
\end{prop}

{\nobf Proof:}
Given $g: \Ss^2 \to \cI_{\bfone}$,
define $\hat g: \Ss^2 \times \Ss^1 \to \Ss^3$
by $\hat g(p,t) = \tilde\fF_{g(p)}(t)$,
where we identify $\Ss^1 = \RR/\ZZ$.
Let $N(g)$ be the degree of $\hat g$.
Clearly $N(g)$ is invariant by homotopy
and $N(g_1 \ast g_2) = N(g_1) + N(g_2)$
(where $\ast$ is the operation defining $\pi_2$).
We prove that $|N(g_0)| = 1$, completing the proof of the proposition.

From Lemma \ref{lemma:nocommontangent} above,
$\tilde\fF_{g(p)}(t) = \bfone$ implies either $t = 0$
or $p = \bfn$ (the north pole) and $t = 1/2$.
For the purpose of computing the degree of $\hat g_0$,
we deform it to define another map $h: \Ss^2 \times \Ss^1 \to \Ss^3$
corresponding to closed curves coinciding with $g_0(p)$
except at a small neighborhood of $t = 0$,
where they cross the $xz$ plane at a point
$(\cos(\epsilon), 0 ,\sin(\epsilon))$, $\epsilon > 0$.
Thus the only preimage of $\bfone$ under $h$ is $(\bfn,1/2)$
and we are left with verifying
that it is topologically transversal.

Alternatively, again from Lemma \ref{lemma:nocommontangent} above,
the preimages under $\hat g_0$ of $-\bfone$
are exactly $(\bfs,1/2)$, $(\bfn,1/4)$ and $(\bfn,3/4)$.
We are left with verifying that all three are topologically transversal
and that the sign of the first is different from the sign
of the last two, again implying $|N(g_0)| = 1$.


Unfortunately, $g_0$ is not differentiable at $\bfs$
but it does admit directional derivatives: that is enough.
We go back to the construction
of $\tilde B_\alpha$ and $\tilde\Gamma_\alpha$ in order
to compute directional derivatives of $g_0$.
For $t \in [-1/12,1/12]$, we may translate the geometric
description above as:
\begin{gather*}
\tilde B_\alpha(t) = \exp\left(\frac{\alpha}{2}\bfj\right)
\exp\left( u(\alpha) t \bfk \right) \exp\left( v(\alpha) \bfj\right), \\
u(\alpha) = 6 \arccos\left( \frac{\cos\alpha}{\sqrt{4-\sin^2 \alpha}} \right),
\quad
v(\alpha) = - \meio \arcsin\left( \meio \sin\alpha \right).
\end{gather*}
Notice that $u'(0) = 0$, $v'(0) = -1/4$.
Let $\tilde B_\alpha^{\bullet}(t)$ be the derivative
of $\tilde B_\alpha(t)$ with respect to $\alpha$;
we have
$(\tilde B_0(t))^{-1} \tilde B_0^{\bullet}(t) = w(t)$
where the auxiliary function $w$ is defined by
\[ w(t) =
\frac{1}{4} \left( \left( 2 \cos(4\pi t) - 1 \right) \bfj +
\left( - 2 \sin(4\pi t) \right) \bfi \right). \]
It is now easy to obtain similar formulas for other intervals
and to deduce that
\[ (\tilde B_0(t))^{-1} \tilde B_0^{\bullet}(t) =
\begin{cases}
w(t-t_0), &
t \in \left [ t_0 - \frac{1}{12}, t_0 + \frac{1}{12} \right ],
t_0 = \frac{k}{3}, \\
-w(t-t_0), &
t \in \left [ t_0 - \frac{1}{12}, t_0 + \frac{1}{12} \right ],
t_0 = \frac{k}{3} + \frac{1}{6}, 
\end{cases}
\quad
k \in \ZZ.
\]
Thus $(\tilde B_0(t))^{-1} \tilde B_0^{\bullet}(t)$,
as a function of $t$,
performs three full turns around the origin
in the plane spanned by $\bfi$ and $\bfj$.
We now have that $(B_0(t;t+\meio))^{-1} B_0^{\bullet}(t;t+\meio)$
performs one full turn around the origin
when $t$ goes from $0$ to $1/3$.

Recall that $\tilde\Gamma_\alpha(t) = \tilde B_\alpha(t) \bfh^{-1}$
and therefore
$\Gamma_\alpha(t;t+\meio) = \bfh B_\alpha(t;t+\meio) \bfh^{-1}$
and we have that
$(\Gamma_0(t;t+\meio))^{-1} \Gamma_0^{\bullet}(t;t+\meio)$
performs one full turn around the origin
when $t$ goes from $0$ to $1/3$,
but now in the plane spanned by $\bfih$ and $\bfj$.
A similar computation shows that, when $t$ goes from $0$ to $1/3$,
$(\Gamma_\pi(t;t+\meio))^{-1} \Gamma_\pi^{\bullet}(t;t+\meio)$
also performs one full turn around the origin
in the same plane.

Translating this back to $g_0$ shows that when $p$
describes a small circle around either the south or the north pole,
$g_0(p)(1/2)$ describes a small simple closed curve
around $(g_0(\bfs))(1/2) = (g_0(\bfn))(1/2) = e_1$,
with $(g_0(p))'(1/2) \approx (g_0(\bfs))'(1/2) = (g_0(\bfn))'(1/2) = e_2$.
The reader should check in Figure \ref{fig:Gamma5}
that this is indeed the case: said simple closed curve is drawn
clockwise when we go left to right along
either the second row from the bottom (around $\bfs$)
or the second from the top (around $\bfn$).
It follows from Lemma \ref{lemma:nocommontangent}
(and can be checked in Figure \ref{fig:Gamma5})
that the same holds for other values of $t$.
Thus, for instance,
when $p$ describes a left to right small circle around the north pole,
$g_0(p)(1/4)$ and $g_0(p)(3/4)$ both describe
simple closed curves around $e_1$,
also oriented clockwise.

Finally, translating these results to $\hat g_0$,
the image under $\hat g_0$ of a small sphere around
$(\bfn,1/2)$ wraps once around $\bfone = \hat g_0((\bfn,1/2))$,
proving topological transversality at this point.
Similarly, the image of a small sphere around
$(\bfs,1/2)$, $(\bfn,1/4)$ or $(\bfn,3/4)$
wraps once around $-\bfone$;
the orientation is different for the first point because,
from the point of view of $\Ss^2$,
a left to right circle near $\bfn$ and
a left to right circle near $\bfs$ have opposive orientations.
\qed

\section{Adding loops}

\label{sect:homotosur}

In this section we present a few facts
related to adding loops to a curve (or family of curves),
including the proof of Proposition \ref{prop:homotosur}.
This is of course similar to the proof
of the Hirsch-Smale theorem (\cite{Hirsch}, \cite{Smale}).
The reader familiar with Gromov's ideas will also recognize
this as an easy instance of the h-principle (\cite{EM}, \cite{Gromov});
others will be reminded of Thurston's method for performing
the sphere eversion by corrugations (\cite{Levy}).

\begin{figure}[ht]
\psfrag{xn}{$\times n$}
\psfrag{t0}{$t_0$}
\psfrag{t0pm}{$t_0 - \epsilon$, $t_0 + \epsilon$}
\psfrag{t0mm}{$t_0 - 2\epsilon$}
\psfrag{t0pp}{$t_0 + 2\epsilon$}
\begin{center}
\epsfig{height=25mm,file=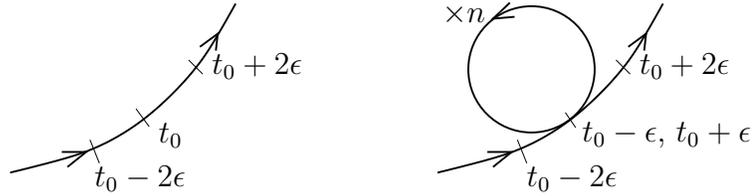}
\end{center}
\caption{Curves $\gamma$ and $\gamma^{[t_0\#n]}$.}
\label{fig:addlooph}
\end{figure}

We need a precise version for the notion
of adding $n$ loops at a point $t_0$ of a curve $\gamma$
as in Figure \ref{fig:addlooph}.
For $\gamma \in \cI$, $t_0 \in (0,1)$ and $n$ a positive integer
let $\gamma^{[t_0\#n]} \in \cI$ be defined by
\[ \gamma^{[t_0\#n]}(t) = 
\begin{cases}
\gamma(t), & 0 \le t \le t_0 - 2\epsilon; \\
\gamma(2t - t_0 + 2\epsilon), & t_0 - 2\epsilon \le t \le t_0 - \epsilon; \\
\fF_\gamma(t_0) \nu_n\left( \frac{t-t_0+\epsilon}{2\epsilon} \right), &
t_0 - \epsilon \le t \le t_0 + \epsilon; \\
\gamma(2t - t_0 - 2\epsilon), & t_0 + \epsilon \le t \le t_0 + 2\epsilon; \\
\gamma(t), & t_0 +2\epsilon \le t \le 1;
\end{cases} \]
here $\epsilon > 0$ is taken sufficiently small.
In words, we follow $\gamma$ normally almost until $t_0$:
we then run a little in order to have room to insert
$\nu_n$ (appropriately moved to the correct position),
run a little again and then continue as $\gamma$.
Since reparametrizations of curves
are not particularly interesting
(the group of orientation-preserving diffeomorphisms of $[0,1]$ is contractible)
the precise value of $\epsilon$ is not particularly interesting either.

For $t_0 = 0$ and $t_0 = 1$, the definition
must be slightly different so that endpoints remain untouched.
Thus, for instance,
\[
\gamma^{[0\#n]}(t) = 
\begin{cases}
\nu_n\left( \frac{t}{\epsilon} \right), &
0 \le t \le \epsilon; \\
\gamma(2t - 2\epsilon), & \epsilon \le t \le 2\epsilon; \\
\gamma(t), & 2\epsilon \le t \le 1.
\end{cases}
\]
Notice that if $\gamma \in \cI_z$ then
$\gamma^{[t_0\#n]} \in \cI_{(-\bfone)^n z}$.
Also, if $\gamma$ is locally convex then so is $\gamma^{[t_0\#n]}$,
with $\tot(\gamma^{[t_0\#n]}) = 2\pi n + \tot(\gamma)$.
We use the notation $\gamma^{[t_0\#n_0; t_1\#n_1]}$ for
$(\gamma^{[t_0\#n_0]})^{[t_1\#n_1]} = (\gamma^{[t_1\#n_1]})^{[t_0\#n_0]}$.
Also,
given $t_0: K \to (0,1)$ and $f: K \to \cI_{Q}$ continuous functions,
let $f^{[t_0\#n]}: K \to \cI_{Q}$ be defined by
$f^{[t_0\#n]}(p) = (f(p))^{[t_0(p)\#n]}$.

Notice that in $\cI_{Q}$
it is easy to introduce a pair of loops at any point of the curve.

\begin{lemma}
\label{lemma:addloopi}
Let $K$ be a compact set, $Q \in SO_3$ and $n$ a positive even integer.
Let $t_0: K \to (0,1)$ and $f: K \to \cI_{Q}$ be continuous functions.
Then $f$ and $f^{[t_0\#n]}$ are homotopic.
\end{lemma}

\begin{figure}[ht]
\begin{center}
\epsfig{height=9mm,file=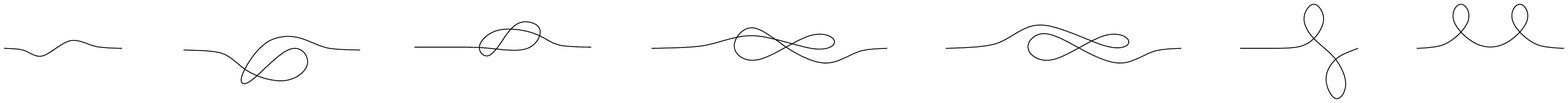}
\end{center}
\caption{How to add two small loops to a curve in $\cI_{Q}$.}
\label{fig:addloopi}
\end{figure}

{\nobf Proof:}
The process is illustrated in Figure \ref{fig:addloopi};
in the final step one of the loops becomes big,
goes around the sphere and shrinks again.
We do not think that an explicit formula is helpful.
\qed

Somewhat similarly, for locally convex curves we have the following.

\begin{lemma}
\label{lemma:addloopl}
Let $K$ be a compact set, $Q \in SO_3$ and $n > 1$ an integer.
Let $t_0: K \to (0,1)$ and $f: K \to \cL_{Q}$ be continuous functions.
Then $f^{[t_0\#n]}$ and $f^{[t_0\#(n+2)]}$ are homotopic,
i.e., there exists $H: [0,1] \times K \to \cL_Q$
with $H(0,p) = f^{[t_0\#n]}(p)$, $H(1,p) = f^{[t_0\#(n+2)]}(p)$.
We may furthermore assume that
$2\pi n + \tot(f(p)) \le \tot(H(s,p)) \le 2\pi (n+2) + \tot(f(p))$.
\end{lemma}

{\nobf Proof:}
We use the path from $\nu_n$ to $\nu_{n+2}$ in $\cL_I$
constructed in Lemma \ref{lemma:2to4}
to modify the curve $f(p)$ in the interval
$[t_0(p) - \epsilon, t_0(p) + \epsilon]$ only.
\qed

The above statement actually holds for $n = 1$
(by borrowing a little elbow room from the surrounding curve)
but we do not need this observation.
On the other hand, the statement critically fails for $n = 0$.
Indeed, we know
(from Fenchel, Little and others; \cite{Fenchel}, \cite{Little})
that $\nu_1$ and $\nu_3$ are not in the same connected
component of $\cL_I$.

We now need a construction corresponding to spreading loops
along the curve, as in Figure \ref{fig:addloop}.
For $\gamma \in \cI$ and $n > 0$, define 
\[ \gamma^{[\#(2n)]} = {\gamma}^{[t_0\#1;t_1\#2;t_2\#2;
\cdots;t_{n-1}\#2;t_n\#1]}; \quad t_j = \frac{j}{n}.  \]
We assume here that the same $\epsilon > 0$ is used for each loop
so that
\[ t \in [0,1], |t - t_j| \le \epsilon \quad \Rightarrow \quad
\gamma^{[\#(2n)]}(t) = \fF_\gamma(t_j)
\nu_1\left(\frac{t-t_j}{\epsilon}\right). \]
We now have $2n$ loops attached to the curve $\gamma$.
We need another construction, however, to smooth out
the remaining small arcs of $\gamma$,
in order to define a curve $\gamma^{[\flat(2n)]}: [0,1] \to \Ss^2$
which, for sufficiently large $n$,
will be similar to $\gamma^{[\#(2n)]}$
and locally convex (even when the original curve $\gamma$ is not).
For $1 \le j \le n$, let 
\[ t_{j,0} = t_{j-1} + \frac{7\epsilon}{8}, \quad
t_{j,\frac12} = \frac{t_{j-1}+t_j}{2} = \frac{2j-1}{2n}, \quad
t_{j,1} = t_j - \frac{7\epsilon}{8} \]
and define
\[ \gamma^{[\flat(2n)]}(t) = \gamma^{[\#(2n)]}(t), \qquad
t \notin \bigcup_{0 < j \le n} (t_{j,0}, t_{j,1}). \]
It remains to define the arcs
$\gamma^{[\flat(2n)]}: [t_{j,0}, t_{j,1}] \to \Ss^2$.
Notice that
$\fF_{\gamma}(t_{j,\frac12}) = \fF_{\gamma^{[\#(2n)]}}(t_{j,\frac12})$
and, for large $n$,
\[ \fF_{\gamma^{[\#(2n)]}}(t_{j,0}) \approx 
\fF_{\gamma}(t_{j,\frac12})
\fF_{\nu_1}\left(-1/8\right), \quad 
\fF_{\gamma^{[\#(2n)]}}(t_{j,1}) \approx 
\fF_{\gamma}(t_{j,\frac12})
\fF_{\nu_1}\left(+1/8\right); \]
We may therefore apply Lemma \ref{lemma:ellipse}
to conclude that there exists a unique arc of ellipse
parametrized by projective arc-length
$\gamma^{[\flat(2n)]}: [t_{j,0}, t_{j,1}] \to \Ss^2$ with
\begin{gather*}
\fF_{\gamma^{[\flat(2n)]}}(t_{j,0}) = \fF_{\gamma^{[\#(2n)]}}(t_{j,0}),
\qquad
\fF_{\gamma^{[\flat(2n)]}}(t_{j,1}) = \fF_{\gamma^{[\#(2n)]}}(t_{j,1}), \\
\gamma^{[\flat(2n)]}(t_{j,\frac12}) = \gamma^{[\#(2n)]}(t_{j,\frac12}).
\end{gather*}
This completes the definition of $\gamma^{[\flat(2n)]}$.
For clarity, we rephrase it more informally.
Draw circles tangent to the curve at the points
$t_j = j/n$, one for $j = 0$ or $j = n$,
two for other values of $j$.
In order to jump from one circle to the next,
draw an arc of ellipse.

Let $f: K \to \cI_z$ be a continuous function.
Define $f^{[\flat(2n)]}: K \to \cI_z$ by
$f^{[\flat(2n)]}(p) = (f(p))^{[\flat(2n)]}$.
Given $f$, for sufficiently large $n$
the function $f^{[\flat(2n)]}$ is well defined and continuous 
and its image is contained in $\cL_z$.

\begin{lemma}
\label{lemma:spread}
Let $K$ be a compact set, $f: K \to \cI_{z}$
and $t_0: K \to (0,1)$ continuous maps.
Then, for sufficiently large $n$,
the following properties hold.
\begin{enumerate}[(a)]
\item{The image of $f^{[\flat(2n)]}$ is contained in $\cL_{z}$.}
\item{The function $f^{[\flat(2n)]}$ is homotopic to  $f$,
i.e., there exists $H_b: [0,1] \times K \to \cI_{z}$
such that $H_b(0,\cdot) = f^{[\flat(2n)]}$ 
and $H_b(1,\cdot) = f$.}
\item{The function $f^{[\flat(2n)]}$ is homotopic to  $f^{[t_0\#(2n)]}$,
i.e., there exists $H_c: [0,1] \times K \to \cI_{z}$
such that $H_c(0,\cdot) = f^{[\flat(2n)]}$ 
and $H_c(1,\cdot) = f^{[t_0\#(2n)]}$.}
\item{If the image of $f$ is contained in $\cL_{z}$ then
the image of the homotopy $H_c$ is also contained in $\cL_{z}$.}
\end{enumerate}
\end{lemma}

Notice that even if the image of $f$ is contained in $\cL_{z}$
we do not claim that $f$ is homotopic to $f^{[\flat(2n)]}$ in $\cL_{z}$;
on the contrary, we shall soon see that this is not always the case.
Notice also that Proposition \ref{prop:homotosur} follows directly
from this lemma.

{\nobf Proof: }
Item (a) follows from the remarks above.

For item (c), notice first that
the functions $f^{[t_0\#(2n)]}$ and $f^{[\#(2n)]}$
are homotopic: the homotopy consists of merely rolling
loops along the curve. More precisely,
for $\tilde t_j(s) = sj/n + (1-s)t_0$,
define
\[ H_1(s,p) = (f(p))^{[\tilde t_0(s)\#1;\tilde t_1(s)\#2;
\cdots;\tilde t_{n-1}(s)\#2;\tilde t_n(s)\#1]}. \]

We next verify that, for sufficiently large $n$,
the functions $f^{[\#(2n)]}$ and $f^{[\flat(2n)]}$
are homotopic. 
Let \(Q_j(p) = (\fF_{f(p)}(t_{j,\frac12}))^{-1} \in SO_3\),
where $t_{j,0}, t_{j,\frac12}, t_{j,1}$ are
as in the construction of $f^{[\flat(2n)]}$.
We have
\begin{gather*}
Q_j(p) \fF_{(f(p))^{[\flat(2n)]}}(t_{j,0}) =
Q_j(p) \fF_{(f(p))^{[\#(2n)]}}(t_{j,0}) \approx 
\Pi(\exp(-\pi\bfkh/8)), \\
Q_j(p) \fF_{(f(p))^{[\flat(2n)]}}(t_{j,1}) =
Q_j(p) \fF_{(f(p))^{[\#(2n)]}}(t_{j,1}) \approx 
\Pi(\exp(+\pi\bfkh/8)).
\end{gather*}
Thus, for sufficiently large $n$, the arcs
\[
Q_j(p) (f(p))^{[\flat(2n)]}, 
Q_j(p) (f(p))^{[\#(2n)]}: [t_{j,0}, t_{j,1}] \to \Ss^2 \]
are \textit{graphs},
in the sense that the first coordinate
$x: [t_{j,0}, t_{j,1}] \to [x_{-}, x_{+}]$
is an increasing diffeomorphism (with $x_{\pm} \approx \pm 1/2$)
and $y$ and $z$ can be considered functions of $x$.
Since the space of increasing diffeomorphisms of an interval
is contractible, we may construct a homotopy from
$f^{[\#(2n)]}$ to a suitable reparametrization
$f_1$ of $f^{[\#(2n)]}$ in each $[t_{j,0}, t_{j,1}]$
for which the function $x$ above is the same as for $f^{[\flat(2n)]}$.
We may then join $f_1$ and $f^{[\flat(2n)]}$
by performing a convex combination followed by projection to $\Ss^2$,
completing the proof of (c).

For item (d), we observe that if the curves are locally convex
then both constructions above remain in the space of locally convex curves.

Finally, for item (b),
we know from Lemma \ref{lemma:addloopi}
that $f$ is homotopic to $f^{[t_0\#(2n)]}$
and from item (c) that $f^{[t_0\#(2n)]}$ is homotopic to $f^{[\flat(2n)]}$.
\qed

As we shall see later, a function $f: K \to \cL_z \subset \cI_z$
may be homotopic to a constant in $\cI_z$
but not in $\cL_z$.
The following proposition shows that this changes if we add loops.

\begin{prop}
\label{prop:easyloop}
Let $n$ be an even positive integer.
Let $K$ be a compact set and 
let $f: K \to \cL_z \subset \cI_z$ a continuous function.
Then $f$ is homotopic to a constant in $\cI_z$
if and only if
$f^{[t_0\#n]}$ is homotopic to a constant in $\cL_z$.
\end{prop}

{\nobf Proof: }
In $\cI_z$, $f$ and $f^{[t_0\#n]}$ are homotopic,
proving one implication.
For the other implication,
let $H: K \times [0,1] \to \cI_z$ be a homotopy
with $H(\cdot,0) = f$, $H(\cdot,1)$ constant.
By Lemma \ref{lemma:spread}, for sufficiently large even $m$,
the image of $H^{[\flat(2m)]}$ is contained in $\cL_z$.
This implies that $f^{[\flat(2m)]}$
is homotopic in $\cL_z$ to a constant.
From Lemma \ref{lemma:spread}, $f^{[t_0\#(2m)]}$
is homotopic to $f^{[\flat(2m)]}$ in $\cL_z$ and therefore
the proposition is proved for large even $n$.
The general case now follows from Lemma \ref{lemma:addloopl}.
\qed

A map $f: K \to \cL_Q$ is \textit{loose}
if $f$ is homotopic to $f^{[t_0\#2]}$ (in $\cL_Q$)
and \textit{tight} otherwise.
Lemma \ref{lemma:addloopl} shows that $f^{[t_0\#2]}$ is loose.
If $K = \{p_0\}$ consists of a single point
then a function $f: K \to \cL_Q$
is essentially a curve $\gamma_0 = f(p_0)$;
$f$ is then tight if and only if $\gamma_0$ is convex
(for $Q=I$ this follows from the results of Little; 
otherwise from Anisov, Shapiro and Shapiro;
\cite{Little}, \cite{Anisov}, \cite{Shapiro2}, \cite{ShapiroM}).
As we shall see, the map $g_0$ constructed above is tight:
this observation will be crucial.

From now on we consider that a main question is,
given $f: K \to \cL_Q$, to decide whether $f$ is loose or tight.
We shall see a few key examples of tight maps
and we shall prove that large classes of maps are loose.
It is sometimes important to have estimates of
the total curvature during the homotopy.

\begin{lemma}
\label{lemma:loose}
Let $f_0, f_1: K \to \cL_Q$ be homotopic with $f_0$ loose.
Then $f_1$ is loose.
\end{lemma}

{\nobf Proof: }
Let $t_0: K \to (0,1)$ be a continuous function
and let $H: [0,1] \times K \to \cL_Q$ be a homotopy
from $f_0$ to $f_1$.
Let $H^{[t_0\#2]}$ be defined by
$H^{[t_0\#2]}(s,p) = (H(s,p))^{[t_0\#2]}$;
clearly, this is a homotopy from
$f_0^{[t_0\#2]}$ to $f_1^{[t_0\#2]}$.
Thus, if $f_0$ is homotopic to $f_0^{[t_0\#2]}$
then $f_1$ is homotopic to $f_1^{[t_0\#2]}$, as desired.
\qed

We finish this section with a more complicated lemma
which allows us to see that many maps $f: K \to \cL_Q$ are loose.

\begin{lemma}
\label{lemma:looppox}
Let $Q \in SO_3$.
Let $K$ be a compact manifold
and $f: K \to \cL_Q$ a continuous map.
Assume that:
\begin{itemize}
\item{ $t_0 \in (0,1)$ and 
$t_1, t_2, \ldots, t_{J}: K \to (0,1)$
are continuous functions with $t_0 < t_1 < t_2 < \cdots < t_{J}$; }
\item{ $K = \bigcup_{1 \le j \le J} U_j$,
where $U_j \subset K$ are open sets; }
\item{ there exist continuous functions $g_j: U_j \to \cL_Q$
such that, for all $p \in U_j$, we have $f(p) = (g_j(p))^{[t_j(p)\#2]}$. }
\end{itemize}
Then $f$ is loose, i.e.,
there exists $H: [0,1] \times K \to \cL_Q$
with $H(0,p) = f(p)$, $H(1,p) = (f(p))^{[t_0\#2]}$.
We may furthermore assume that
\[ \tot(f(p)) \le \tot(H(s,p)) \le 4\pi + \tot(f(p)). \]
\end{lemma}

{\nobf Proof: }
Our proof proceeds by induction on $J$.
For $J = 1$ we have $U_1 = K$ and therefore $f = g_1^{[t_1\#2]}$,
which is known from Lemma \ref{lemma:addloopl} to be loose.
The estimate on the total curvature is also given in the lemma;
notice that sliding a loop between $t_0$ and $t_1$
does not affect total curvature.

Assume now that $J > 1$.
Let $W \subset U_J$ be an open set
whose closure is contained in $U_J$ and such that
$K = W \cup \bigcup_{1 \le j \le J-1} U_j$.
We now slide the loop in $t_J$ to position $t_{J-1}$ in $W$,
allowing for the loop to stop elsewhere for $p \in U_J \smallsetminus W$.
More precisely, let $u: K \to [0,1]$
be a continuous function with
$u(p) = 1$ for $p \in W$ and
$u(p) = 0$ for $p \notin U_J$.
Define $H_J: [0,1] \times K \to \cL_Q$ by
\[ H_J(s,p) = \begin{cases}
f(p),& p \notin U_J, \\
g_J(p)^{[((1-u(p)s)t_J(p) + u(p)s t_{J-1}(p))\#2]}, & p \in U_J. 
\end{cases} \]
Notice that $\tot(H_J(s,p)) = \tot(f(p))$.
Let $\hat f(p) = H_J(1,p)$,
$\hat U_j = U_j$ for $j < J-1$ and $\hat U_{J-1} = U_{J-1} \cup W$;
the hypotheses of the Lemma apply to $\hat f$ with a smaller value
of $J$ and therefore $\hat f$ is loose.
By Lemma \ref{lemma:loose}, so is $f$.
The estimate on total curvature also follows.
\qed

\section{Multiconvex curves}

\label{sect:onion}

A curve $\gamma \in \cL_{z}$ is \textit{multiconvex}
of multiplicity $k$ if
there exist $0 = t_0 < t_1 < \cdots < t_k = 1$ such that
\begin{enumerate}[(a)]
\item{$\fF_\gamma(t_i) = I$ for $i < k$;}
\item{the restrictions $\gamma|_{[t_{i-1},t_i]}$
are convex arcs for $1 \le i \le k$.}
\end{enumerate}
Notice that for $i < k$ these restrictions are then simple closed curves
(see Figure \ref{fig:multiconvex}).
Let $\cM_{k} \subset \cL_{z}$ be
the set of multiconvex curves of multiplicity $k$.

Notice that $\nu_k$ is multiconvex of multiplicity $k$.
A curve $\gamma \in \cL_z$ is multiconvex of multiplicity $1$
if and only if it is convex,
so that for $z = -\bfone$ we have $\cM_1 = \cL_{-\bfone,c}$.
By Lemma \ref{lemma:total}, if $\gamma \in \cM_k$
then $2(k-1)\pi < \tot(\gamma) < 4k\pi$.
It is easy to see that, for $k$ odd,
$\cM_{k} \ne \emptyset$ if and only if $z$ is convex;
similarly, for $k$ even,
$\cM_{k} \ne \emptyset$ if and only if $-z$ is convex.

In \cite{S1}, other submanifolds $\cF_k \subset \cL$
(of \textit{flowers} of order $k$, or of $2k-1$ petals) are introduced
which play a role somewhat similar to $\cM_k$.
For results up to this point, it is indeed largely a matter of taste
to use multiconvex curves or flowers.
For the final part of the paper, however,
multiconvex curves work better.

\begin{lemma}
\label{lemma:Mkcontract}
Let $z \in \Ss^3$.
Let $k$ be a positive integer.
The closed subset $\cM_{k} \subset \cL_{z}$
(if non-empty)
is a contractible submanifold of codimension $2k-2$
with trivial normal bundle.
\end{lemma}

{\nobf Proof:}
Assume that $-z$ is convex.
Consider the geodesic $\rho \subset \Ss^2$ passing through
$\pm e_1$ and $\pm e_3$ (where $e_1 = (1,0,0)$).
After a projective transformation
we may assume that
any convex curve $\gamma \in \cL_{-z}$
crosses $\rho$ transversally once for some $t \in (0,1)$.

We first define open sets $U_k \subset \cL_{(-\bfone)^k z}$.
A curve $\gamma \in \cL_{(-\bfone)^k z}$ belongs to $U_k$ if and only if:
\begin{enumerate}[(a)]
\item{all intersections between $\gamma$ and $\rho$ are transversal;}
\item{there are exactly $2k$ values
\[ 0 = t_0 < t_{\meio} < t_1 < \cdots < t_{k-1} < t_{k - \meio} < 1 \]
of $t \in [0,1)$ for which $\gamma(t) \in \rho$;}
\item{consecutive intersections $\gamma(t_j)$ and $\gamma(t_{j\pm \meio})$
are distinct;}
\item{arcs of $\gamma$ between $t_j$ and $t_{j+\meio}$,
$t_{j+\meio}$ and $t_{j+1}$ or
$t_{k - \meio}$ and $1$ are convex.}
\end{enumerate}
Notice that $U_k$ is indeed open and $t_j: U_k \to [0,1]$
are continuous functions.
We may continuously define the arguments $\theta_j$
of the points $\gamma(t_j(\gamma))$ by
\[
\gamma(t_j(\gamma)) = \cos(\theta_j(\gamma)) e_1 +
\sin(\theta_j(\gamma)) e_3, \quad
\theta_0(\gamma) = 0, \]
and, for integer $j$,
\[ \theta_j < \theta_{j+\meio} < \theta_j + \pi, \quad
\theta_{j+\meio} > \theta_{j+1} > \theta_{j+\meio} - \pi. \]
Also, $\langle \gamma'(t_j(\gamma)), e_2 \rangle$
is positive if $j$ is an integer (and negative otherwise).
For $j$ an integer, we continuously define the argument $\eta_j$
of $\gamma'(t_j(\gamma))$ by
\begin{gather*}
\gamma'(t_j(\gamma)) = r \left(
\cos(\eta_j(\gamma)) e_2 + \sin(\eta_j(\gamma)) n \right), \\
n = - \sin(\theta_j(\gamma)) e_1 + \cos(\theta_j(\gamma)) e_3, \quad
-\pi/2 < \eta_j(\gamma) < \pi/2, \quad r > 0. 
\end{gather*}

Now define $M_k: U_k \to \RR^{2k-2}$ by
\[ M_k(\gamma) = ( \theta_1(\gamma), \eta_1(\gamma), 
\theta_2(\gamma), \eta_2(\gamma), \ldots,
\theta_{k-1}(\gamma), \eta_{k-1}(\gamma) ). \]
The smooth map $M_k$ is a submersion and
$\cM_k$ is the inverse image of $0 \in \RR^{2k-2}$.
This proves that $\cM_k$ is a smooth submanifold
of codimension $2k-2$ and trivializes its normal bundle.

Finally, we must prove that $\cM_k$ is contractible.
The fact that $\cL_{-z,c} = \cM_1(z)$ is contractible is well known
(\cite{Anisov} and \cite{ShapiroM}, Lemma 5).
The subset $\hat\cM_k$ of $\cM_k$ of curves for which $t_j = j/k$
is naturally identified with $(\cM_1(-\bfone))^{(k-1)} \times (\cM_1(z))$
(reparametrize $\gamma|_{[t_j,t_{j+1}]}$ to define $\gamma_j \in \cM_1$)
and therefore is also contractible.
But $\hat\cM_k \subset \cM_k$ is a deformation retract:
just use piecewise affine functions to reparametrize each curve
so that $t_j = j/k$ (for all $j$).
\qed



The previous result allows us to use each $\cM_k$
to define an element $m_{2k-2} \in H^{2k-2}(\cL_z;\ZZ)$
by counting intersections with $\cM_k$.
For de Rham cohomology, for instance,
we consider Thom's form in the (trivial) normal bundle
to $\cM_k$ and use the identification of this bundle
with a tubular neighborhood of $\cM_k$ to define
a closed $(2k-2)$-form $\omega$ which
is a representative of $m_{2k-2}$.
Thus, if $f: K \to \cL_z$ is a smooth map from an oriented
compact $(2k-2)$-dimensional manifold $K$
to $\cL_z$ which is transversal to $\cM_k$
then the integral of the pull back of $\omega$
over $K$ equals the number of intersections
of $f$ with $\cM_k$, counted with sign.
If the map is not smooth or not (topologically) transversal
we may perturb it so that it becomes both smooth and transversal:
the number of intersections is still well defined.
We denote this integer by $m_{2k-2}(f)$.
The elements $m_{2k-2} \in H^{2k-2}(\cL_z)$
will turn out to be the ``extra'' cohomology (as compared to $\cI_z$);
compare with Corollary \ref{coro:oldtheo}.

As we shall see in Lemma \ref{lemma:g2I} below,
the map $g_0$ (introduced in Section \ref{sect:2to4})
is tight and satisfies $m_2(g_0) = \pm 1$ but $m_2(g_0^{[t_0\#2]}) = 0$.
For $N \in H^2(\cL_{+\bfone};\ZZ)$ as defined
in Proposition \ref{prop:g0generator}, we have
$N(g_0) = N(g_0^{[t_0\#2]}) = \pm 1$
(from Proposition \ref{prop:g0generator} and Lemma \ref{lemma:addloopi}).
It follows that $m_2$ and $N$ span a copy of $\ZZ^2$
in $H^2(\cL_{+\bfone};\ZZ)$
and that $g_0$ and $g_0^{[t_0\#2]}$ span a copy of $\ZZ^2$
in $\pi_2(\cL_{+\bfone}) = H_2(\cL_{+\bfone};\ZZ)$.
Compare this with Corollary \ref{coro:oldtheo}:
we shall later see that $m_2$ and $N$ actually span
$H^2(\cL_{+\bfone};\ZZ)$ and that
$g_0$ and $g_0^{[t_0\#2]}$ actually span $\pi_2(\cL_{+\bfone})$.

In the following lemma the sphere $\Ss^2$ will
be the compact manifold (usually called $K$)
in the domain of a map.
Let $\bfs = -e_3$, $\bfn = e_3$ be
the south and north pole, respectively.
The base point of $\Ss^2$ is $\bfs$.

\begin{lemma}
\label{lemma:g2I}
There exist maps $g_s: \Ss^2 \to \cL_{\bfone}$,
$s \in [0, \meio)$, such that:
\begin{enumerate}[(a)]
\item{$g_s(\bfs) = \nu_{2}$ and $g_s(\bfn)$
is a reparametrization of $\nu_{4}$;}
\item{if $g_s(p)$ is multiconvex then $p = \bfs$ or $p = \bfn$;}
\item{$g_s$ is topologically transversal to $\cM_2$ at $\bfs$;}
\item{if $t \in [1-s,1]$ then $(g_s(p))(t) = \nu_2(t)$;}
\item{if $\tilde\fF_{g_s(p)}(t) \in \cC_{\bfkh}$, $t \in (0,1)$,
then either $p = \bfs$, $p = \bfn$ or $t \in [1-s,1]$;}
\item{given $t \in (0,1-s)$, the map $p \mapsto \tilde\fF_{g_s(p)}(t)$
is topologically transversal to $\cC_{\bfkh}$;}
\item{the maps $g_s$ are all homotopic to $g_0$;}
\item{the maps $g_s$ satisfy $m_2(g_s) = \pm 1$ and are all tight.}
\end{enumerate}
\end{lemma}

We reiterate that
the map $g_0$ is the same one introduced in Section \ref{sect:2to4}.
The maps $g_s$ for $s > 0$ may be informally described
as modifications of $g_0$ by forcing $g_s$ to coincide
with $\nu_2$ for $t \in [1-s,1]$ (as in item (d)).
The proof of the lemma thus splits into three parts:
we first prove that $g_0$ satisfies all the desired properties,
we next construct $g_s$ for $s > 0$ and we finally verify
that the properties remain true for $s > 0$.
More precisely, we first construct $g_s$ for small positive $s$
by making small adjustments to $g_0$;
for large $s$ we use projective transformations.
The constructions of $g_s$ are rather explicit,
and it should be noted that many arbitrary choices 
are made during the construction, sometimes in order to facilitate
some later argument.
The maps $g_s$ will be the building blocks in the construction
of the maps $h_{2k-2}$ in Lemma \ref{lemma:h2k-2} below.

The sign ambiguity in the last item comes from the fact
that we were not too careful to define either
a standard transversal orientation to $\cM_2$ or
a standard orientation for $K = \Ss^2$, the domain of $g_s$.
Recall that $\cC_{\bfkh} \subset \Ss^3$ is a subgroup
defined before Lemma \ref{lemma:nocommontangent}.

{\nobf Proof:}
We first prove that $g_0$ satisfies the desired properties.
Items (a) and (d) are immediate and item (e) follows directly
from Lemma \ref{lemma:nocommontangent}.
Item (b) is a direct consequence of item (e).
For items (c) and (f) we first notice
that we are talking about isolated points.
Indeed, from (b), $p = \bfs$ is the only point 
for which $g_0(p) \in \cM_2$.
Similarly, from (e), $p = \bfs$ is the only point
for which $\fF_{g_0(p)}(t) \in \cC_{\bfkh}$.
In either case we need to study what happens
when $p$ goes around $\bfs$, drawing a small circle.
For (c), we need to prove that $g_0(p)$ will go around $\cM_2$ once,
or, equivalently, that $M_2(g_0(p))$ will go around the origin once.
For (f), we need to prove that $\fF_{g_0(p)}(t)$ will go around
the circle $\cC_{\bfkh}$ once.
Notice that both observations are rather clear from Figure \ref{fig:Gamma5}.

Go back to the topological transversality
argument in Proposition \ref{prop:g0generator}.
Recall that $(B_0(t;t+\meio))^{-1} B_0^{\bullet}(t;t+\meio)$
performs one full turn around the origin
when $t$ goes from $0$ to $1/3$.

Let $p = (\cos\theta \sin\alpha, \sin\theta \sin\alpha, -\cos\alpha)$;
let $M_2$ and $U_2$ be as in Lemma \ref{lemma:Mkcontract}.
For sufficiently small $\alpha$, we have $g_0(p) \in U_2$.
From the above computations, for sufficiently small $\alpha$,
$M_2(g_0(p))$ also performs a full turn around the origin
when $\theta$ goes from $0$ to $2\pi$, completing the proof of (c).
Item (f) follows similarly from these computations for $t = 1/2$;
we notice that, by continuity,
the number of turns of $\tilde\fF_{g_0(p)}(t)$ around $\cC_{\bfkh}$
must be constant; this completes the proof of item (f).

Notice now that items (a) through (f) imply that
$g_0$ intersects $\cM_2$ topologically transversally and exactly once.
Thus, for $m_2 \in H^2$ as above, $m_2(g_0) = \pm 1$.
On the other hand, $g_0^{[t_0\#2]}$ can not possibly
intersect $\cM_2$ and therefore $m_2(g_0^{[t_0\#2]}) = 0 \ne m_2(g_0)$.
Thus $g_0$ and $g_0^{[t_0\#2]}$ are not homotopic and $g_0$ is therefore tight.
This completes the proof of item (h), of the case $s = 0$
and of the remarks preceding the statement of the lemma.



We now construct $g_s$ for $s > 0$, $s$ small.
By compactness, there exists $\epsilon_1 > 0$ such that
for all $p \in \Ss^2$,
the arc $g_0(p)|_{[-2\epsilon_1,+2\epsilon_1]}$ is convex.
Here we interpret $g_0(p)$ as a $1$-periodic function
from $\RR$ to $\Ss^2$.
Recall that $\tilde\fF_{g_0(p)}(0) = \tilde\fF_{\nu_2}(0) = \bfone$.
Again by compactness, there exists $\epsilon_2 \in (0,\epsilon_1/4)$
such that, for all $t \in [0,\epsilon_2]$ and for all $p \in \Ss^2$,
we have 
\[ \left(\fF_{\nu_2}(t)\right)^{-1} \fF_{g_0(p)}(\epsilon_1),
\left(\fF_{g_0(p)}(-\epsilon_1)\right)^{-1} \fF_{\nu_2}(-t)
\in \Bruhat_{(13);2}. \]
We want to define $\hat g: \Ss^2 \to \cL_{\bfone}$ with
\[ (\hat g(p))(t) = \begin{cases}
\nu_2(t), & t \in [0,\epsilon_2] \cup [1-\epsilon_2, 1], \\
(g_0(p))(t), & t \in [\epsilon_1, 1-\epsilon_1].
\end{cases}
\]
As in Lemma \ref{lemma:ellipse}, for each $p \in \Ss^2$,
there exist ellipses $\cE_+$ and $\cE_-$ and parametrizations
by projective arc-length
$\gamma_+: [\epsilon_2, \epsilon_1] \to \cE_+ \subset \Ss^2$
and
$\gamma_-: [1-\epsilon_1, 1-\epsilon_2] \to \cE_- \subset \Ss^2$
such that
\begin{gather*}
\gamma_+(\epsilon_2) = \nu_2(\epsilon_2), \quad
\gamma'_+(\epsilon_2) = \nu_2'(\epsilon_2), \quad
\fF_{\gamma_+}(\epsilon_1) = \fF_{g(p)}(\epsilon_1), \\
\gamma_-(1-\epsilon_2) = \nu_2(1-\epsilon_2), \quad
\gamma'_-(1-\epsilon_2) = \nu_2'(1-\epsilon_2), \quad 
\fF_{\gamma_-}(1-\epsilon_1) = \fF_{g(p)}(1-\epsilon_1), 
\end{gather*}
and, furthermore, $\cE_+$ and $\cE_-$ osculate
the circle $\cC_{\bfkh}$ at
$\nu_2(\epsilon_2)$ and $\nu_2(1-\epsilon_2)$, respectively.
The ellipses and parametrizations are uniquely and continuously defined.
Complete the definition of $\hat g$ by
\[ (\hat g(p))(t) = \begin{cases}
\gamma_+(t), & t \in [\epsilon_2, \epsilon_1], \\
\gamma_-(t), & t \in [1-\epsilon_1, 1-\epsilon_2].
\end{cases}
\]
Notice that since $5$ points define a conic,
the ellipses $\cE_\pm$ have no tangency point to $\cC_{\bfkh}$
besides $\nu_2(\pm\epsilon_2)$.

For $s \le 2\epsilon_2$, set
$g_{s}(p)(t) = (\fF_{\nu_2}(\epsilon_2))^{-1} \hat g(t-\epsilon_2)$.
We claim that the function $g_{s}$ has all the required properties.
Item (a) is obvious.
Item (d) holds by construction
and item (e) follows from the last observation
in the previous paragraph;
item (b) now follows.
Topological transversality (items (c) and (f))
is handled as for $s = 0$.
Item (g) follows either from an explicit computation
or from the contractibility of $\cL_{z,c}$
and item (h) now follows.
This completes the proof of the theorem for $s \le {2\epsilon_2}$.

For $c \in \RR$, let
\[ A(c) = \begin{pmatrix} 1 & c & c^2/2 \\
0 & 1 & c \\ 0 & 0 & 1 \end{pmatrix} =
\exp\left( c \begin{pmatrix} 0 & 1 & 0 \\
0 & 0 & 1 \\ 0 & 0 & 0 \end{pmatrix} \right). \] 
Notice that $\pi(A(c)) e_1 = e_1$, $\pi(A(c)) I = I$ and
$\pi(A(c))$ takes the circle $\cC_{\bfkh}$ to itself.
Let $\hat g_c: \Ss^2 \to \cL_{\bfone}$
be given by $\hat g_c(p) = \pi(A(c)) \circ g_{\epsilon_2}(p)$.
Given $c$, each arc $(\hat g_c(p))_{[1-\epsilon_2,1]}$
is a fixed parametrization of an arc of $\cC_{\bfkh}$.
Changing $c$, that arc can be taken to have any required length.
In other words, given $s \in (\epsilon_2,\meio)$
there exists a unique $c \in \RR$
for which $\fF_{\hat g_c(p)}(1-\epsilon_2) = \fF_{\nu_2}(1-s)$:
define $g_s$ by suitably reparametrizing this $\hat g_c$.
All the required properties follow by construction.
\qed


Let $z \in \Ss^3$ with $-z$ convex.
We are now ready to construct tight maps
$h_{2k-2}: \Ss^{2k-2} \to \cL_{(-\bfone)^k z}$
corresponding to the spheres attached to $\cI_{\Pi(z)}$
to obtain $\cL_{\Pi(z)}$
as in Theorems \ref{theo:main} and \ref{theo:mainplus}.
Alternatively, $h_{2k-2}$ define the ``extra'' generators
of $H_\ast(\cL_{\Pi(z)};\ZZ)$ (compared to $H_\ast(\cI_{\Pi(z)};\ZZ)$);
see Corollary \ref{coro:oldtheo}.

\begin{lemma}
\label{lemma:h2k-2}
Let $k > 1$ be a positive integer.
Let $z \in \Ss^3$ with $-z$ convex.
There exist tight maps $h_{2k-2}: \Ss^{2k-2} \to \cL_{(-\bfone)^k z}$ which:
\begin{enumerate}[(a)]
\item{intersect $\cM_{k}$ exactly once and topologically transversally; }
\item{do not intersect $\cM_{k'}$ for $k' \ne k$; }
\item{satisfy $m_{2k-2}(h_{2k-2}) = \pm 1$; }
\item{are homotopic to a constant as maps
$\Ss^{2k-2} \to \cI_{(-\bfone)^k z}$. }
\end{enumerate}
\end{lemma}

{\nobf Proof:}
Let $\DD^2 \subset \RR^2$ be the closed disk of radius $1$.
We first construct a function
$\hat h: (\DD^2)^{(k-1)} \to \cL_{(-\bfone)^k z}$.
Consider $\epsilon_0 > 0$ such that if $t \in (0, \epsilon_0)$ then
$- (\tilde\fF_{\nu_1}(kt))^{-1} z$ is convex.
Let $s_0 = k(1+\epsilon_0)$ and
$z_0 = \tilde\fF_{\nu_{s_0}}(\frac{1}{k}) = - \tilde\fF_{\nu_1}(\epsilon_0)$.
We shall have $(\hat h(0))(t) = \nu_{s_0}(t)$
for $t \in [0,\frac{k-1}{k}]$
and $\tilde\fF_{\hat h(p)}(\frac{i}{k}) = z_0^i$
for all $p$ and for all integers $i < k$
(in other words, we are starting the definition here;
the expression ``we shall have'' is to be understood as:
``here is yet another property of $\hat h$,
which is clearly consistent with what we demanded before'').

Let $\gamma_k: [\frac{k-1}{k},1] \to \Ss^2$ be a convex arc
with $\tilde\fF_{\gamma_k}(\frac{k-1}{k}) = z_0^{k-1}$,
$\tilde\fF_{\gamma_k}(1) = (-\bfone)^k z$.
We shall have $(\hat h(p))(t) = \gamma_k(t)$
for all $p \in (\DD^2)^{(k-1)}$ and $t \in [\frac{k-1}{k},1]$.
Let $s_1 = 1-\frac{s_0}{2k} = \frac{1 - \epsilon_0}{2}$
and let $g_{s_1}$ be as in Lemma \ref{lemma:g2I}
so that $\tilde\fF_{g_{s_1}(p)}({1-s_1}) = z_0$ for all $p \in \Ss^2$.
Recall that $\tilde\nu_4 = g_{s_1}(\bfn)$ is a reparametrization of $\nu_4$
with $\tilde\nu_4(t) = \nu_2(t)$ for all $t > 1-s_1$.
Define $w: \DD^2 \to \Ss^2$ by 
\[ w(r \cos \theta, r \sin \theta) = \begin{cases}
(\cos\theta \sin(4r), \sin\theta \sin(4r), -\cos(4r)),
& r \le \pi/4, \\
(0,0,1), & r \ge \pi/4.
\end{cases} \]
Consider $p = (p_1,p_2, \ldots, p_{k-1}) \in (\DD^2)^{(k-1)}$,
$p_i \in \DD^2$.
For $t \in [\frac{i-1}{k}, \frac{i}{k}]$ 
let $t_i = (1-s_1)k \left( t - \frac{i-1}{k} \right)$;
if $|p_i| \le \frac{\pi}{4}$
we shall have $(\hat h(p))(t) = \Pi(z_0^{i-1}) g_{s_1}(w(p_i))(t_i)$.
Let $\tilde\nu_{8k}$ be a reparametrization of $\nu_{8k}$ with 
$\tilde\nu_{8k}(t) = \nu_2(t)$ for all $t > 1-s_1$.
Let $\tilde g: [\frac{\pi}{4},\frac{7}{8}] \to \cL_{+\bfone}$ be a path
from $\tilde g(\frac{\pi}{4}) = \tilde\nu_4$
to $\tilde g(\frac{7}{8}) = \tilde\nu_{8k}$ satisfying
$\tilde g(\tau)(t) = \nu_2(t)$ for all $t > 1-s_1$;
notice that $\tilde\fF_{\tilde g(\tau)}(1-s_1) = z_0$
for all $\tau \in [\frac{\pi}{4},\frac{7}{8}]$.
For $t \in [\frac{i-1}{k}, \frac{i}{k}]$ and $t_i$ as above
we shall have
\[ (\hat h(p))(t) = \begin{cases}
\Pi(z_0^{i-1}) \tilde g(|p_i|)(t_i), & |p_i| \in [\frac{\pi}{4},\frac{7}{8}], \\
\Pi(z_0^{i-1}) \tilde \nu_{8k}(t_i), & |p_i| \ge \frac{7}{8},
\end{cases} \]
completing the construction of $\hat h$.

From Lemma \ref{lemma:g2I},
$\hat h(0) \in \cM_{k}$,
$\hat h(p) \in \cM_{k}$ implies $p = 0$
and $\hat h$ is topologically transversal to $\cM_k$,
with a single intersection at $p = 0$.
Consider
\[ \hat h_1 = \hat h|_{\partial((\DD^2)^{(k-1)})}:
\partial((\DD^2)^{(k-1)}) \to \cL_{(-\bfone)^k z}. \]
We apply Lemma \ref{lemma:looppox} to prove that $\hat h_1$ is loose:
here $K = \partial((\DD^2)^{(k-1)})$
(which is homeomorphic to $\Ss^{(2k-3)}$),
$J = k-1$, $t_j = \frac{j}{k} - \frac{1}{2k}$
and $(p_1,p_2, \ldots, p_{k-1}) \in U_j$ if
$i \ne j$ implies $|p_i| > \frac{7}{8}$.
There exists therefore a homotopy
$H: [0,1] \times \partial((\DD^2)^{(k-1)}) \to \cL_{(-\bfone)^k z}$
with $H(0,p) = \hat h_1(p)$ and $H(1,p) = (\hat h_1(p))^{[t_k\#2]}$,
$t_k = 1 - \frac{1}{2k}$.

The homotopy $H$ may be assumed to be disjoint from $\cM_k$.
In order to see this we give estimates on the total curvature.
The total curvature of $\tilde\nu_{8k}$ equals $16k\pi$,
and the total curvature of its restriction to $[0,1-s_1]$
is greater than $(16k-2)\pi$.
If $p \in \partial((\DD^2)^{(k-1)})$ we have
at least one index $j$ for which $|p_j| = 1$;
for such $j$, the total curvature in the interval
$[\frac{j-1}{k},\frac{j}{k}]$ is greater than $(16k-2)\pi$.
The total curvature of $\hat h(p)$ is therefore greater
than $(16k-2)\pi$.
By Lemma \ref{lemma:looppox}, we may construct $H$ as above
with $\tot(H(s,p)) > (16k-2)\pi > 4k\pi$
(for all $s \in [0,1]$ and all $p \in \partial((\DD^2)^{(k-1)})$)
and therefore $H$ is disjoint from $\cM_k$, as claimed.

Let $\DD^{(2k-2)}$ be the closed disk of dimension $2k-2$ and radius $1$.
There exists a homeomorphism from $\DD^{(2k-2)}$ to 
$(\{0\} \times (\DD^2)^{(k-1)}) \cup ([0,1] \times \partial((\DD^2)^{(k-1)}))$.
Compose this homeomorphism with $\hat h$ and $H$
to define a map $\tilde h: \DD^{(2k-2)} \to \cL_{(-\bfone)^k z}$
with $\tilde h(0) \in \cM_{k}$;
$\tilde h(p) \in \cM_{k}$ implies $p = 0$;
$\tilde h$ is topologically transversal to $\cM_k$,
with a single intersection at $p = 0$;
$\tot(\tilde h(p)) > 16k\pi$ for all $p \in \Ss^{(2k-3)}$.
We furthermore have $\tilde h(p) = \gamma^{[t_k\#2]}$
for all $p \in \Ss^{(2k-3)}$ (for some $\gamma$).
More precisely, let $Q_k = \fF_{\gamma_k}(t_k)$;
after a reparametrization we may assume that,
for all $p \in \Ss^{(2k-3)}$,
$(\tilde h(p))(t_k + \epsilon_1 \tau)  = Q_k\nu_2(\tau+\meio)$
for $|\tau| \le \meio$,
where $\epsilon_1 \in (0,\frac{1}{4k})$ is a small positive constant.

The sphere $\Ss^{(2k-2)}$ is homeomorphic to
$S = (\{0,1\} \times \DD^{(2k-2)}) \cup ([0,1] \times \Ss^{(2k-3)})$.
Let $\gamma_2: [0,1] \to \cL_{\bfone}$ with
$\gamma_2(0) = \nu_2$, $\gamma_2(1) = \nu_4$.
Define $\tilde h_{2k-2}: S \to \cL_{(-\bfone)^k z}$ by 
\[ \tilde h_{2k-2}(s,p)(t) = \begin{cases}
(\tilde h(p))(t), &
t \notin (t_k - \frac{\epsilon_1}{2}, t_k + \frac{\epsilon_1}{2}), \\
(\tilde h(p))(t), & s = 0, \\
(\tilde h(p))^{[t_k\#2]}(t), & s = 1, \\
Q_k \gamma_2(s)(\tau+\meio), & |p| = 1,\; t = t_k + \epsilon_1 \tau,\;
\tau \in [-\meio,\meio].
\end{cases} \]
We claim that $\tilde h_{2k-2}$ is homotopic to a point
in $\cI_{(-\bfone)^k z}$, or, equivalently (Proposition \ref{prop:easyloop}),
that $\tilde h_{2k-2}^{[t_k\#2]}$ can be extended to a map from
$[0,1] \times \DD^{(2k-2)}$ to $\cL_{(-\bfone)^k z}$.
Indeed, after a reparametrization,
$\tilde h^{[t_k\#2]}: \DD^{(2k-2)} \to \cL_{(-\bfone)^k z}$
may be assumed to satisfy
$(\tilde h^{[t_k\#2]}(p))(t_k + \epsilon_1 \tau)
 = \tilde Q(p) \nu_2(\tau+\meio)$
for all $p \in \DD^{(2k-2)}$,
where $\tau \in [-\meio,\meio]$
and $\tilde Q: \DD^{(2k-2)} \to SO_3$ satisfies
$\tilde Q(p) = Q_k$ if $|p| = 1$.
Define $\bar h: [0,1] \times \DD^{(2k-2)} \to \cL_{(-\bfone)^k z}$
by
\[ \bar h(s,p)(t) = \begin{cases}
(\tilde h^{[t_k\#2]}(p))(t), &
t \notin (t_k - \frac{\epsilon_1}{2}, t_k + \frac{\epsilon_1}{2}), \\
\tilde Q(p) \gamma_2(s)(\tau+\meio), &
t = t_k + \epsilon_1 \tau,\; \tau \in [-\meio,\meio];
\end{cases} \]
up to reparametrization, $\bar h$ is the desired extention.

Identifying $S$ with $\Ss^{2k-2}$,
the function $\tilde h_{2k-2}: \Ss^{2k-2} \to \cL_{(-\bfone)^k z}$
thus satisfies item (d).
By construction, the only multiconvex curves
in its image are $\nu_{k} = \tilde h_{2k-2}(0,0)$
and $\tilde h_{2k-2}(1,0)$, which is a reparametrization of $\nu_{k+2}$.
Define $h_{2k-2}$ by perturbing $\tilde h_{2k-2}$
near $(1,0)$ so as to avoid $\cM_{k+2}$;
by transversality, this can be done:
the codimension of $\cM_{k+2}$ is larger than the dimension of $\Ss^{2k-2}$.
Item (b) is therefore satisfied.
Topological transversality in item (a) also follows by construction
and by items (c) and (f) of Lemma \ref{lemma:g2I}.
Finally, $m_{2k-2}(h_{2k-2}) = \pm 1$ follows from items
(a) and (b), proving item (c).
\qed

The following corollary sums up some of the topological differences
between the spaces $\cL_z$ and $\cI_z$ which we have proved in this section.

\begin{coro}
\label{coro:coho}
Consider $z \in \Ss^3$ with $-z$ convex.
For $k \ge 1$,
the elements 
\( m_{2k-2} \in H^{2k-2}(\cL_{(-\bfone)^kz}) \)
do not belong to the image of
\( i^\ast = H^\ast(i):
H^\ast(\cI_{(-\bfone)^kz}) \to H^\ast(\cL_{(-\bfone)^kz})\).
The maps $h_{2k-2}: \Ss^{2k-2} \to \cL_{(-\bfone)^kz}$ define 
non-zero elements in the kernel of
\( i_\ast = \pi_{2k-2}(i): \pi_{2k-2}(\cL_{(-\bfone)^kz}) \to
\pi_{2k-2}(\cI_{(-\bfone)^kz}) \).
\end{coro}

The aim of the rest of the paper is to prove that these are,
in a sense, the \textit{only} differences.

Our final lemma in this section is an easy consequence
of the previous results and will be used later.

\begin{lemma}
\label{lemma:Mkgohome}
Let $B$ be a compact manifold of dimension $n+1$
with boundary $\partial B = K$.
Let $f_0: B \to \cL_z$ be a continuous map
with $f_0|_{K}$ disjoint from all $\cM_k$.
Then there exists a continuous map $f_1: B \to \cL_z$
with $f_0|_{K} = f_1|_{K}$ and $f_1$ disjoint from all $\cM_k$.
\end{lemma}

{\nobf Proof: }
For each $k$, let $V_k \subset \cL_z$ be
a closed tubular neighborhood of $\cM_k$.
We may assume the sets $V_k$ to be disjoint
and $f_0|_{K}$ to be disjoint from the sets $V_k$.
We may furthermore assume the map $h_{2k-2}: \Ss^{2k-2} \to \cL_z$
to be topologically transversal to $\partial V_k$
so that $D_k = h_{2k-2}^{-1}(V_k) \subset \Ss^{2k-2}$
is a disk around the south pole.
Since $\cM_k$ is a contractible Hilbert manifold,
$\cM_k$ is homeomorphic to the Hilbert space $\cH$.
Let $\DD^{2k-2} \subset \RR^{2k-2}$ be the closed unit ball of radius $1$.
Let $\psi_k = (\psi_{k,1}, \psi_{k,2}): V_k \to \cH \times \DD^{2k-2}$
be a homeomorphism taking $\cM_k$ to $\cH \times \{0\}$;
we may assume that $\psi_{k,1} \circ h_{2k-2}$ is constant equal to $0$
and that $\phi_k = \psi_{k,2} \circ h_{2k-2}$ is 
a homeomorphism between $D_k$ and $\DD^{2k-2}$.
For $c \in \RR$ and $\gamma \in V_k$,
let $\mu(c,\gamma) \in V_k$ be such that
$\psi_{k,1}(\mu(c,\gamma)) = c \psi_{k,1}(\gamma)$ and
$\psi_{k,2}(\mu(c,\gamma)) = \psi_{k,2}(\gamma)$.

Let $\tilde D_k \subset D_k$ be the inverse image
under $\phi_k$ of the closed disk of radius $1/2$.
Let $\rho: \tilde D_k \to \Ss^{2k-2}$ be a continuous map
coinciding with the identity on $\partial\tilde D_k$
and avoiding the south pole.
Define 
\[ f_1(p) = \begin{cases}
f_0(p), & f_0(p) \notin \bigcup_k V_k, \\
\mu(2|\psi_{k,2}(f_0(p))| - 1, f_0(p)),
& f_0(p) \in V_k, |\psi_{k,2}(f_0(p))| \ge \meio, \\
h_{2k-2}(\rho(\phi_k^{-1}(\psi_{k,2}(f_0(p))))),
& f_0(p) \in V_k, |\psi_{k,2}(f_0(p))| \le \meio.
\end{cases} \]
The map $f_1$ satisfies the required conditions.
\qed

\section{Grafting}

\label{sect:graft}

In this section we introduce the process of \textit{grafting} curves;
similar ideas are considered in \cite{SZ} and \cite{Z}.

Subintervals $[t_0, t_1] \subset (0,1)$
will be called \textit{arcs}:
an arc $[t_0, t_1]$ is \textit{graftable}
for a locally convex curve $\gamma$ if there exists
a projective transformation taking $\gamma$ to $\gamma_1$
and real numbers $\theta_0, \theta_1, \phi_0, \phi_1$
with $-\pi/2 < \phi_0 < 0 < \phi_1 < \pi/2$ and
\begin{equation}
\label{eq:graft}
\tilde\fF_{\gamma_1}(t_0) = e^{\theta_0\bfk/2} e^{\phi_0\bfj/2}, \quad
\tilde\fF_{\gamma_1}(t_1) = e^{\theta_1\bfk/2} e^{\phi_1\bfj/2} e^{\pi\bfi/2}.
\end{equation}
Recall that $\Pi(e^{\theta\bfk/2})$ is a rotation by $\theta$
around the $z$ axis: 
\[ \Pi(e^{\theta\bfk/2}) =
\begin{pmatrix}
\cos\theta  & -\sin\theta & 0 \\
\sin\theta  & \cos\theta & 0 \\
0 & 0 & 1
\end{pmatrix};
\]
similarly, $\Pi(e^{\phi\bfj/2})$ is a rotation by $\phi$ around the $y$ axis.

Translating into a more geometric language,
equation \ref{eq:graft} says that
$\gamma_1$ is tangent at $t_i$ to the circle $z = -\sin(\phi_i)$;
in both cases the orientation of $\gamma_1$ is consistent
with a locally convex orientation of the circles.

Notice that the existence of the desired projective transformation
depends only on the value of $\fF_\gamma(t_0;t_1)$.
A matrix $Q \in SO_3$ is \textit{graftable}
if $\fF_\gamma(t_0;t_1) = Q$ implies that
the arc $[t_0,t_1]$ is graftable for $\gamma$.

\begin{lemma}
\label{lemma:graftbruhat}
Let $Q \in SO_3$:
$Q$ is graftable if and only if
$Q$ belongs to one of the Bruhat cells below:
\[ \Bruhat_{(13);1}, \Bruhat_{(13);4}, \Bruhat_{(13);7},
\Bruhat_{(123);3}, \Bruhat_{(123);5},
\Bruhat_{(132);5}, \Bruhat_{(132);6}, \Bruhat_{e;5}. \]
\end{lemma}

{\nobf Proof: }
The definition is clearly invariant under
projective transformation and therefore if $Q_1$ and $Q_2$
belong to the same Bruhat cell
then $Q_1$ is graftable if and only if $Q_2$ is.

We need therefore to check which Bruhat cells are touched
by 
\[ \fF_{\gamma_1}(t_0;t_1) =
e^{-\phi_0\bfj/2}
e^{(\theta_1-\theta_0)\bfk/2} e^{\phi_1\bfj/2} e^{\pi\bfi/2}. \]
A straightforward computation (or a few figures)
show that if $\theta_0 = \theta_1$ then
$\fF_{\gamma_1}(t_0;t_1) =
e^{(\phi_1-\phi_0)\bfj/2} e^{\pi\bfi/2} \in \Bruhat_{(13);7}$.
We may keep $\theta_0$ fixed and change $\theta_1$
and the Bruhat cell will cycle.
If $\phi_0 + \phi_1 < 0$ the cell will go through
$P_{(13);7}$, $P_{(123);5}$, $P_{(13);1}$ and $P_{(123);3}$
(and back);
notice that the cells $\Bruhat_{(13);7}$ and $\Bruhat_{(13);1}$
are open and correspond to intervals while the cells
$\Bruhat_{(123);5}$ and $\Bruhat_{(123);3}$
have dimension $2$ and correspond to transition points.
If $\phi_0 + \phi_1 > 0$ the cell will go through
$P_{(13);7}$, $P_{(132);6}$, $P_{(13);4}$ and $P_{(132);5}$.
In the special case $\phi_0 + \phi_1 = 0$
the only transition point is $P_{e;5}$.
\qed

We give a few examples which will be used again later:
\begin{gather*}
Q_{0,1} = \Pi(e^{\pi\bfj/6}), \quad
Q_{1,1} = \Pi(e^{\pi\bfk/2}e^{-\pi\bfj/{12}}e^{\pi\bfi/2}), \\
Q_{0,4} = \Pi(e^{\pi\bfj/{12}}), \quad
Q_{1,4} = \Pi(e^{\pi\bfk/2}e^{-\pi\bfj/{6}}e^{\pi\bfi/2}), \\
Q_{0,7} = \Pi(e^{\pi\bfj/8}), \quad
Q_{1,7} = \Pi(e^{-\pi\bfj/8}e^{\pi\bfi/2}).
\end{gather*}
We clearly have that $Q_{0,\ell}^{-1}Q_{1,\ell} \in \Bruhat_{(13);\ell}$;
also, if $\fF_{\gamma_1}(t_i) = Q_{i,\ell}$ then
$\gamma_1$ satisfies the conditions in equation \ref{eq:graft}.

A minor difficulty is that the choices of
$\theta_0, \theta_1, \phi_0, \phi_1$
and of the projective transformation
should be uniform.
The following lemma addresses this issue.

\begin{lemma}
\label{lemma:graftbruhatu}
Let $\ell$ be equal to $1$, $4$ or $7$.
Let $Q_0, Q_1 \in SO_3$ be such that
\[ Q_0^{-1} Q_1 \in \Bruhat_{(13);\ell}. \]
Then there exists a unique matrix $U \in \Upum$
such that
$\Pi(Q_{0,\ell}UQ_0^{-1})(Q_i) = Q_{i,\ell}$
for $i = 0, 1$.
\end{lemma}

{\nobf Proof: }
We have
$\Pi(Q_0^{-1})(Q_0) = I$ and 
$\Pi(Q_0^{-1})(Q_1) = Q_0^{-1} Q_1 \in \Bruhat_{(13);\ell}$.
There exists a unique $U \in \Upum$ with
$\Pi(U)(Q_0^{-1} Q_1) = Q_{0,\ell}^{-1}Q_{1,\ell}$
and the result follows.
\qed

Given a curve $\gamma \in \cL$ and
a graftable arc $[t_0,t_1]$ such that
$\fF_\gamma(t_0;t_1) \in \Bruhat_{(13);\ell}$
(where $\ell$ equals $1$, $4$ or $7$)
we shall assume that
$\gamma_1 = \pi(Q_{0,\ell}UQ_0^{-1}) \circ \gamma$
where $U \in \Upum$ is as in Lemma \ref{lemma:graftbruhatu}.

\begin{figure}[ht]
\begin{center}
\epsfig{height=55mm,file=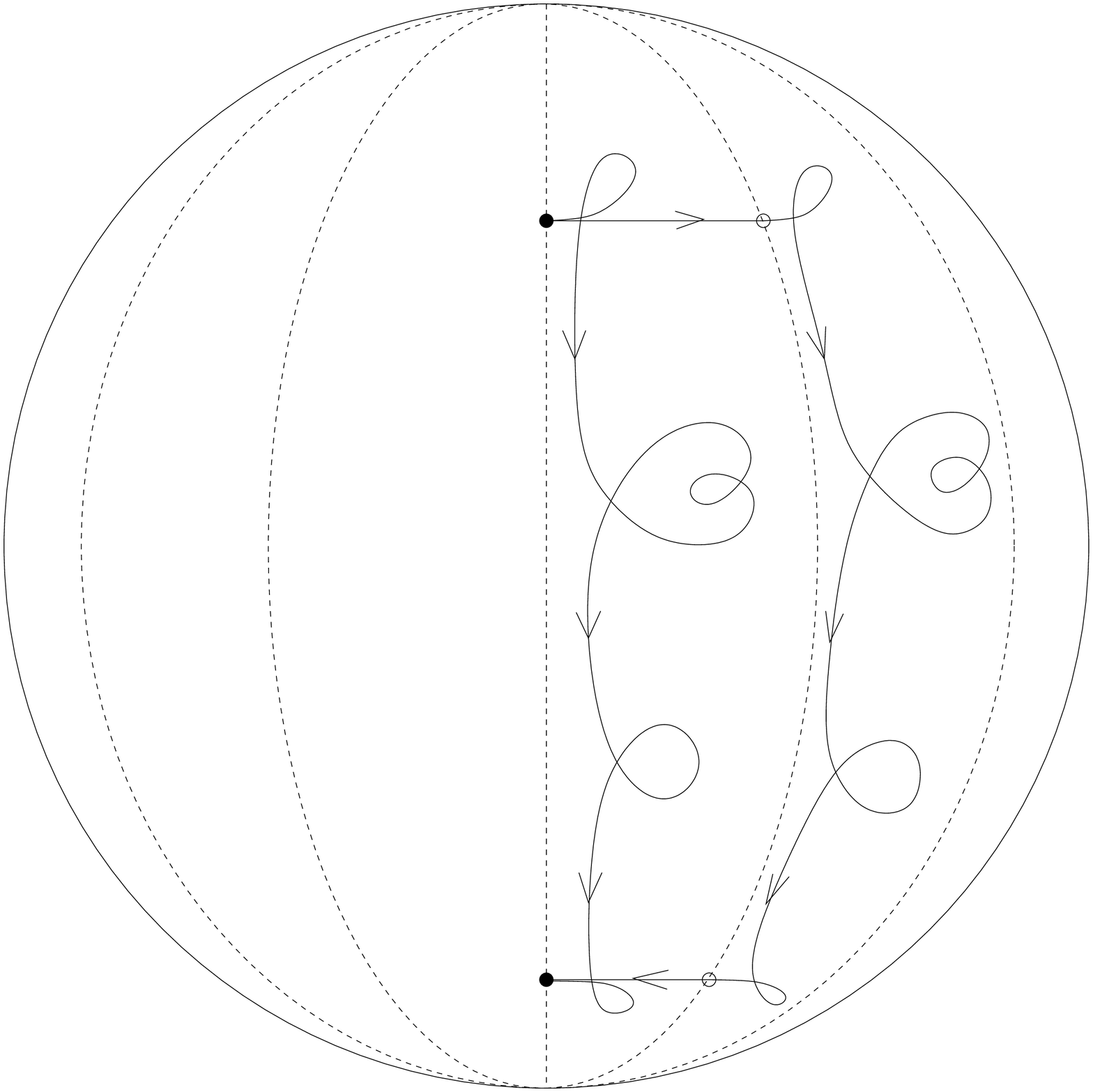}
\end{center}
\caption{Grafting a curve}
\label{fig:loopatendpoints}
\end{figure}

Consider a curve $\gamma_1 \in \cL$ and
a graftable arc $[t_0,t_1]$ such that
equation \ref{eq:graft} is satisfied.
Given a real number $s \ge 0$
we show how to perform a \textit{graft} on the curve $\gamma_1$
around the arc $[t_0,t_1]$
to obtain a curve $\gamma_1^{[(t_0,t_1)\#s]}$.
As usual, write $\tilde\Gamma_1 = \tilde\fF_{\gamma_1}$.
Let $\epsilon > 0$ be a small number.
Define
$\tilde\Gamma_1^{[(t_0,t_1)\#s]}$ by
$\tilde\Gamma_1^{[(t_0,t_1)\#s]}(t) = \tilde\Gamma_1(t)$
for $t \le t_0$ and for $t \ge t_1$; otherwise
\[ \tilde\Gamma_1^{[(t_0,t_1)\#s]}(t) = \begin{cases}
\exp\left(\left(\frac{\theta_0}{2}+\frac{\pi (t-t_0)}{\epsilon}\right)
\bfk \right) e^{\phi_0\bfj/2}, &
t_0 \le t \le t_0+s\epsilon, \\
e^{\pi s\bfk} \tilde\Gamma_1(2t-t_0-2s\epsilon), &
t_0 + s\epsilon \le t \le t_0+2s\epsilon, \\
e^{\pi s\bfk} \tilde\Gamma_1(t), &
t_0 + 2s\epsilon \le t \le t_1-2s\epsilon, \\
e^{\pi s\bfk} \tilde\Gamma_1(2t-t_1+2s\epsilon), &
t_1 - 2s\epsilon \le t \le t_1-s\epsilon, \\
\exp\left(\left(\frac{\theta_1}{2}+\frac{\pi (t_1-t)}{\epsilon}\right)
\bfk \right) e^{\phi_1\bfj/2} e^{\pi\bfi/2}, &
t_1 - s\epsilon \le t \le t_1. \\
\end{cases} \]
Equation \ref{eq:graft} guarantees compatibility and continuity.
Define (of course)
\[ \gamma_1^{[(t_0,t_1)\#s]}(t) =
\Pi(\tilde\Gamma_1^{[(t_0,t_1)\#s]}(t))e_1. \]
A geometric description of this construction is in order:
the curve $\gamma_1$ is cut at $t_0$ and $t_1$,
the arc from $t_0$ to $t_1$ is rotated by an angle of $2\pi s$
around the $z$ axis
and finally arcs of circle (parallel to the plane $z=0$)
are grafted into the resulting gap (see Figure \ref{fig:loopatendpoints}).
Notice that $\tilde\Gamma_1^{[(t_0,t_1)\#s]}(t)$ is continuous
as a function of $s$ and $t$;
in other words, if $\gamma_1 \in \cL_Q$ then
$s \mapsto \gamma_1^{[(t_0,t_1)\#s]}$ 
is a continuous path from $[0, s_{\max}]$ to $\cL_Q$.

We now define grafting on $\gamma$ provided
$\fF_{\gamma}(t_0;t_1) \in \Bruhat_{(13);\ell}$
where $\ell$ equals $1$, $4$ or $7$.
Let $U$ be as in Lemma \ref{lemma:graftbruhatu} and
$M = Q_{0,\ell}UQ_0^{-1} \in SL_3$ so that
the projective transformation $\pi(M)$ takes
$\gamma$ to $\gamma_1 = \pi(M) \circ \gamma$ satisfying
$\fF_{\gamma_1}(t_i) = Q_{i,\ell}$.
Next graft $\gamma_1$ as above to obtain $\gamma_1^{[(t_0,t_1)\#s]}$
and define $\gamma^{[(t_0,t_1)\#s]} =
\pi(M^{-1}) \circ \gamma_1^{[(t_0,t_1)\#s]}$.
Define a group homomorphism $A: \RR \to SL_3$ by
$A(s) = M^{-1} \Pi(e^{s\pi\bfk}) M = \exp(sa)$, $a \in sl_3$;
notice that $A(1) = \exp(a) = I$.
We have
\[ \gamma^{[(t_0,t_1)\#s]}(t) =
\begin{cases}
\gamma(t), & t \le t_0, \\
\pi(\exp(sa)) \gamma(t), & t_0 + 2s\epsilon \le t \le t_1 - 2s\epsilon, \\
\gamma(t), & t \ge t_1.
\end{cases} \]

If $s$ is a positive integer,
the curves $\gamma^{[(t_0,t_1)\#s]}$ (grafting)
and $\gamma^{[t_0\#s;t_1\#s]}$ (adding loops,
as in Section \ref{sect:homotosur}, just before Lemma \ref{lemma:addloopi})
are very similar.
The only significant difference is that
the ``loops'' in the first curve are closed convex curves
which are not quite circles.
Since the space of closed convex curves
with a given base point is contractible
these loops can easily be made round.




We write $A \Subset B$ if the closure of $A$
is contained in the interior of $B$.
Notice that if $A$ is open and closed then $A \Subset A$.

\begin{lemma}
\label{lemma:B147}
Let $\ell$ equal $1$, $4$ or $7$.
Let $K$ be a compact manifold and $f: K \to \cL$ be a continuous map.
Let $K_0 \subset K$ be a compact subset.
Assume that there exist continuous functions
$t_0 < t_1: K_0 \to (0,1)$ such that, for all $p \in K_0$,
$\fF_{f(p)}(t_0; t_1) \in \Bruhat_{(13);\ell}$.
Let $W_0 \Subset K_0$ be an open set.
Then, for sufficiently small $\epsilon > 0$,
there exist a homotopy $H: [0,1] \times K \to \cL$
and a function $A: [0,1] \times K \to SL_3$ with the following properties:
\begin{enumerate}[(a)]
\item{$t_0(p) + 8\epsilon < t_1(p) - 8\epsilon$ for all $p \in K_0$;}
\item{$A(p) = I$ for all $p \in (K \smallsetminus K_0) \cup W_0$;}
\item{$H(0,p) = f(p)$ for all $p \in K$;}
\item{$H(s,p) = f(p)$ for all $p \in K \smallsetminus K_0$
and all $s \in [0,1]$;}
\item{$H(s,p)(t) = f(p)(t)$ for all $p \in K_0$
and for all $t \notin (t_0(p), t_1(p))$ and all $s \in [0,1]$;}
\item{$H(s,p)(t) = \pi(A(s,p))(f(p)(t))$ for all $p \in K_0$
and for all $t \in (t_0(p) + 8\epsilon, t_1(p) - 8\epsilon)$
and all $s \in [0,1]$;}
\item{$H(1,p) = f(p)^{[(t_0(p)+4\epsilon)\#2;(t_1(p)-4\epsilon)\#2]}$
for all $p \in W_0$.}
\end{enumerate}
\end{lemma}


{\nobf Proof: }
Let $W_1 \Subset K_0$ be an open set such that $W_0 \Subset W_1$.
Let $u: K \to [0,1]$ be a smooth function
with $u(p) = 0$ for $p \notin K_0$ and $u(p) = 1$ for $p \in W_1$.

For $s \in [0,1/2]$, define $H(s,p)$ by grafting $f(p)$:
\[ H(s,p) = f(p)^{[(t_0(p),t_1(p))\#(4su(p))]}. \]
Notice that for $p \in W_1$ we have
$H(1/2,p) = f(p)^{[(t_0(p),t_1(p))\#2]}$.
For $s \in [1/2,1]$ and
$p \notin W_1$ we define $H(s,p) = H(1/2,p)$.
For $p \in W_1$ we use the interval to round
up the loops introduced by grafting.
The margin $W_1 \smallsetminus W_0$ is needed for compatibility
but for $p \in W_0$ we may assume that
$H(1,p) = f(p)^{[(t_0(p)+4\epsilon)\#2;(t_1(p)-4\epsilon)\#2]}$,
completing our proof.
\qed

We are now ready to prove the easier cases of our main theorem;
these are also proved in \cite{SaSha} using different ideas.

\begin{coro}
\label{coro:easytheo}
Let $z \in \Ss^3$ with $\Pi(z) \in
\Bruhat_{(13);1} \cup \Bruhat_{(13);4} \cup \Bruhat_{(13);7}$.
Then the inclusion $\cL_z \subset \cI_z$ is a weak homotopy equivalence.
\end{coro}

Since these spaces are Hilbert manifolds,
they are actually diffeomorphic (\cite{BH}).

{\nobf Proof: }
Let $\ell$ be such that $\Pi(z) \in \Bruhat_{(13);\ell}$.
Let $K$ be a compact manifold and $f: K \to \cL_z$
be a continuous map.
For sufficiently small $\epsilon > 0$ we have
$\fF_{f(p)}(\epsilon, 1-\epsilon) \in \Bruhat_{(13);\ell}$.
Apply Lemma \ref{lemma:B147} to $f$
with $t_0 = \epsilon$, $t_1 = 1-\epsilon$,
$K_0 = W_0 = K$ to deduce that
$f$ is homotopic to $f^{[t_a\#2;t_b\#2]}$ and therefore loose.
It now follows from Proposition \ref{prop:easyloop}
that $f$ is homotopic to a constant in $\cL_z$
if and only if $f$ is homotopic to a constant in $\cI_z$.
Together with Lemma \ref{lemma:spread}, 
this completes the proof.
\qed


Part of Little's Theorem is that the set $\cL_{-\bfone,c}$
of simple locally convex curves is a connected component of $\cL_{-\bfone}$:
Fenchel proves that simple closed locally convex curves
are convex (\cite{Fenchel}; see also \cite{Little} and \cite{ShapiroM}); 
we often use this fact.
The other part of Little's Theorem is that,
once convex curves have been removed,
the sets $\cL_{+\bfone}$ and $\cL_{-\bfone,n}$ are path connected.
This is again a corollary of Lemma \ref{lemma:B147}.

\begin{coro}
\label{coro:easylittle}
The sets $\cL_{+\bfone}$ and $\cL_{-\bfone,n}$ are path connected.
\end{coro}

\begin{figure}[ht]
\begin{center}
\epsfig{height=25mm,file=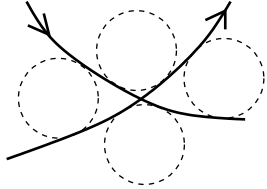}
\end{center}
\caption{A transversal self-intersection}
\label{fig:doublepoint}
\end{figure}

{\nobf Proof: }
Consider a map from $\Ss^0$ (two points)
to either of these spaces (i.e., two curves):
the map is homotopic to a constant in $\cI_{\pm\bfone}$.
Each of the two curves $\gamma_0$ and $\gamma_1$
may be assumed generic and therefore to have
a transversal self-intersection.
As in figure \ref{fig:doublepoint},
near the self-intersection there exist $t_0$ and $t_1$ such that
$\fF_\gamma(t_0; t_1) \in \Bruhat_{(13);7}$
(there are other pairs for which this expression
belongs to any of the other three open cells).
\qed

\section{Good and bad steps}

\label{sect:nextstep}

Given a locally convex curve $\gamma: [t_0,t_1] \to \Ss^2$,
we now define the \textit{next step} function
$\Next_\gamma: [t_0,t_1^{-}] \to [t_0,t_1]$,
or, to make it shorter when $\gamma$ is clear from the context,
$t^{+} = \Next_\gamma(t)$
(the number $t_1^{-} \in [t_0,t_1)$ will also be defined).
Given $t \in [t_0,t_1]$ let $t^{+}$ be the smallest $\tilde t > t$ such that
$\Gamma(t;\tilde t) \notin \Bruhat_{(13),2}$;
if no such $\tilde t$ exists then $t^+$ is undefined.
Since $\Bruhat_{(13),2}$ is an open set,
$\Next_\gamma$ is a continuous function (where defined).
Also, the function $\Next_\gamma$ is strictly increasing
with continuous inverse,
which will be denoted by $t^{-} = \Next_\gamma^{-1}(t)$.
Unless $\gamma: [t_0, t_1] \to \Ss^2$ is convex,
$\Next_\gamma$ is a strictly increasing continuous homeomorphism
from $[t_0, t_1^{-}]$ to $[t_0^{+}, t_1]$.
On the other hand, $\Next_\gamma$ is usually not differentiable
(even if $\gamma$ is smooth).

Geometrically speaking, $t_0^{+}$ is the point
at which the arc $\gamma|_{[t_0,t_0^{+}]}$ is still convex
but about to somehow lose convexity.
This can occur in five different ways
corresponding to five Bruhat cells
to which $\Gamma(t_0; t_0^{+})$ may belong.
The two generic cases are when $\gamma$
is about to leave the hemisphere defined by $\Gamma(t_0)$
(but not at the point $\gamma(t_0)$)
or, conversely, when the geodesic defined by $\Gamma(t)$
passes through $\gamma(t_0)$
(but not aligned with $\gamma'(t_0)$)
so that $\gamma$ is about to enter its own convex hull:
these correspond to $P_{(123);6}$ and $P_{(132);0}$,
in this order, and to the first two diagrams in Figure \ref{fig:5bruhat}:
Notice that these matrices have two inversions
and therefore their Bruhat cells have dimension $2$.
Two more exceptional cases correspond to the matrices
$P_{(23);2}$ and $P_{(12);4}$
which, having one inversion,
correspond to Bruhat cells of dimension $1$.
The two cases correspond to the third and fourth diagram
in Figure \ref{fig:5bruhat}:
the curve may self-intersect by coming back to $\gamma(t_0)$
(but with non-aligned tangent vectors)
or it may tangentially touch the geodesic defined by $\Gamma(t_0)$
(but not at $\gamma(t_0)$).
The fifth and most exceptional case corresponds to $P_{e;0} = I$,
with Bruhat cell of dimension $0$:
the curve may come back to $\gamma(t_0)$
with tangent vectors also aligned
(as in the fifth diagram).

We are particularly interested in this fifth and most degenerate case.
A \textit{step} (for $\gamma$) is an interval $[t_0,t_1]$
with $t_1 = t_0^{+}$.
A step is \textit{bad} if $\Gamma(t_1) = \Gamma(t_0)$
(this is the fifth case) and \textit{good} otherwise.
Notice that the set of good steps is open (in the space of steps).
A curve $\gamma$ is multiconvex if and only if,
for $t_0 = 0$, the steps $[t_0,t_1], [t_1, t_2], \ldots$
are all bad (here $t_j = \Next_\gamma^j(t_0)$).
Conversely, if a curve is complicated
(i.e., not multiconvex) then, again with $t_j = \Next_\gamma^j(t_0)$,
there exists a good step $[t_j,t_{j+1}] \subset [0,1]$

\begin{figure}[ht]
\begin{center}
\epsfig{height=20mm,file=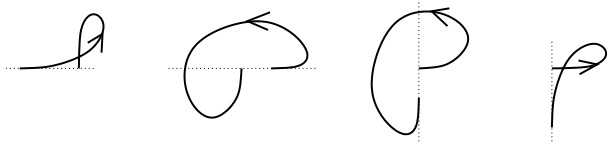}
\end{center}
\caption{Good arcs}
\label{fig:trip}
\end{figure}

A \textit{good arc} for a locally convex curve $\gamma$
is an interval $[t_0, t_1] \subset (0,1)$ such that:
\begin{enumerate}[(a)]
\item{$t_0 < t_1^{-} < t_0^{+} < t_1$;}
\item{if $t \in [t_0, t_1^{-}]$ then
$[t,t^{+}]$ is a good step;}
\item{if $t_0 \le t_a < t_a^{+} < t_b \le t_1$ then
$\fF_\gamma(t_a;t_b)$ is in one
of the following Bruhat cells:
$\Bruhat_{(13);1}$, $\Bruhat_{(13);4}$, $\Bruhat_{(13);7}$,
$\Bruhat_{(123);3}$, $\Bruhat_{(123);5}$,
$\Bruhat_{(132);5}$ or $\Bruhat_{(132);6}$.}
\end{enumerate}
The first three matrices in item (c) correspond of course
to open cells as in Figure \ref{fig:genericbruhat} above;
the last four are shown in Figure \ref{fig:trip}.

A \textit{good pair of arcs} for a locally convex curve $\gamma$
consists of two good arcs
$[t_0, t_1] \subset [\tilde t_0, \tilde t_1] \subset (0,1)$ such that
$\tilde t_0 < t_0 < t_1 < \tilde t_1$ and
if $\fF_\gamma(t_0;t_1) \notin \Bruhat_{(13),1} \cup \Bruhat_{(13),4}$
then $\fF_\gamma(\tilde t_0; \tilde t_1) \in \Bruhat_{(13),7}$.

\begin{lemma}
\label{lemma:goodarc}
Let $K$ be a compact manifold;
let $f: K \to \cL$ be a family of locally convex curves.
Let $t_0, t_1: K \to (0,1)$ be continuous functions
with $\Next_{f(p)}(t_0(p)) = t_1(p)$
such that $[t_0(p), t_1(p)]$ is always a good step.
\begin{enumerate}[(a)]
\item{There exists $\epsilon > 0$ such that
$[t_0(p) - \epsilon, t_1(p) + \epsilon] \subset (0,1)$
is always a good arc.}
\item{For any $\epsilon > 0$ there exist continuous functions
$\tilde t_0 < \hat t_0 < \hat t_1 < \tilde t_1: K \to (0,1)$
with $t_0 - \epsilon < \tilde t_0 < \hat t_0 < t_0$
and $t_1 < \hat t_1 < \tilde t_1 < t_1 + \epsilon$
such that, for all $p \in K$,
$[\hat t_0(p), \hat t_1(p)] \subset [\tilde t_0(p), \tilde t_1(p)]$
is a good pair of arcs for $f(p)$.}
\end{enumerate}
\end{lemma}

{\nobf Proof: }
Consider the union of all allowed cells
in the definition of a good arc
for $\fF_\gamma(t_a;t_b)$ if $t_a < t_a^{+} < t_b$:
\begin{align*}
A_1 &= \Bruhat_{(13);1} \cup \Bruhat_{(13);4} \cup
\Bruhat_{(13);7} \cup \\
&\phantom{=} \cup \Bruhat_{(123);3} \cup \Bruhat_{(123);5} \cup
\Bruhat_{(132);5} \cup \Bruhat_{(132);6} \subset SO_3. 
\end{align*}
Consider also the set of all allowed cells for
$\fF_\gamma(t_a;t_b)$ if $t_a < t_b$:
\[
A_2 = A_1 \cup \Bruhat_{(13);2} \cup
\Bruhat_{(123);6} \cup \Bruhat_{(132);0} \cup
\Bruhat_{(23);2} \cup \Bruhat_{(12);4} \subset SO_3.
\]
We claim that the sets $A_1$ and $A_2$ are both open.
For $A_1$, each of the four $2$-cells included
is sandwiched between two of the included $3$-cells:
$\Bruhat_{(123);3}$ between $\Bruhat_{(13);1}$ and $\Bruhat_{(13);7}$;
$\Bruhat_{(123);5}$ between $\Bruhat_{(13);1}$ and $\Bruhat_{(13);7}$;
$\Bruhat_{(132);5}$ between $\Bruhat_{(13);4}$ and $\Bruhat_{(13);7}$;
$\Bruhat_{(132);6}$ between $\Bruhat_{(13);4}$ and $\Bruhat_{(13);7}$.
For $A_2$, since all four $3$-dimensional cells are included
we need only check that each of the two $1$-dimensional cells included
is surrounded by $2$- and $3$-dimensional cells which are also included:
$\Bruhat_{(23);1}$ is surrounded by
the four $3$-dimensional cells plus
$\Bruhat_{(132);0}$, $\Bruhat_{(123);6}$, 
$\Bruhat_{(123);3}$ and $\Bruhat_{(132);6}$; 
$\Bruhat_{(12);4}$ is surrounded by
the four $3$-dimensional cells plus
$\Bruhat_{(123);6}$, $\Bruhat_{(132);0}$, 
$\Bruhat_{(132);5}$ and $\Bruhat_{(123);5}$. 
This completes the proof of the claim

By hypothesis, if $t_0(p) \le t_a < t_b \le t_1(p)$
then $\fF_{f(p)}(t_a;t_b) \in A_2$.
Thus, using the compactness of $K$ and the fact that $A_2$ is open,
for sufficiently small $\epsilon_2$,
if $t_0(p) -\epsilon_2 < t_a < t_b < t_1(p) + \epsilon_2$
then $\fF_{f(p)}(t_a;t_b) \in A_2$.
In particular, $[t_a,t_a^{+}]$ is a good step.
For any good step $[t_a,t_a^{+}]$ for a locally convex curve $\gamma$
there exists $\epsilon_a > 0$ such that
if $t_a^{+} < t_b < t_a^{+} + \epsilon_a$
then $\fF_\gamma(t_a;t_b) \in A_1$.
Again by compactness, there exists therefore $\epsilon \in (0,\epsilon_2)$
such that, for all $p$,
if $t_0(p) -\epsilon < t_a < t_a^{+} < t_b < t_1(p) + \epsilon$
then $\fF_{f(p)}(t_a;t_b) \in A_1$,
proving (a).

For item (b),
assume arcs parametrized by a constant multiple 
(depending on $p$ only) of arc length.
We claim that for sufficiently small $\epsilon_b > 0$,
we may take
$\tilde t_0(p) = t_0(p) - 2\epsilon_b$,
$\hat t_0(p) = t_0(p) - \epsilon_b$,
$\hat t_1(p) = t_1(p) + \epsilon_b$,
$\tilde t_1(p) = t_1(p) + 2\epsilon_b$:
$[\hat t_0(p), \hat t_1(p)] \subset [\tilde t_0(p), \tilde t_1(p)]$
is a good pair of arcs for $f(p)$.

Let $B = \Bruhat_{(123);6} \cup \Bruhat_{(132);0}
\cup \Bruhat_{(23);2} \cup \Bruhat_{(12);4} \subset A_2 \subset SO_3$;
the set $B$ is a topological manifold of dimension $2$
homeomorphic to $\Ss^1 \times (0,1)$;
the subsets $\Bruhat_{(23);2}, \Bruhat_{(12);4} \subset B$
are closed with disjoint neighborhoods $B_2, B_4 \subset B$,
respectively.
We may assume the closures of $B_2$ and $B_4$ in $B$ to be disjoint.
Let $s: K \to B$ be defined by $s(p) = \fF_{f(p)}(t_0(p);t_1(p))$;
let $U_2 = s^{-1}(B_2)$, $U_4 = s^{-1}(B_4)$.
Define open sets
$U_6 \Subset s^{-1}(\Bruhat_{(123);6})$
and $U_0 \Subset s^{-1}(\Bruhat_{(132);0})$,
so that the sets $U_i$ form an open cover of $K$.
For sufficiently small $\epsilon > 0$,
if $t_0(p) -\epsilon < t_a < t_a^{+} < t_b < t_1(p) + \epsilon$
and $p \in U_6$ (resp. $p \in U_0$)
then $\fF_{f(p)}(t_a;t_b) \in \Bruhat_{(13);4}$
(resp. $\fF_{f(p)}(t_a;t_b) \in \Bruhat_{(13);1}$).
For $p \in U_6 \cup U_0$, therefore,
pairs of arcs will be good.
We must focus on $p \in U_2$ and $p \in U_4$.

Assume $p \in U_2$.
For small $\epsilon > 0$, we may assume that the arc
$[t_0(p) - \epsilon, t_1(p) + \epsilon]$ has at most one self intersection.
Let $V_2 \subset U_2$ be the open set of points $p \in U_2 \subset K$
for which there exist $t_c, t_d \in (t_0(p) - \epsilon, t_1(p) + \epsilon)$
with $t_c < t_d$, $f(p)(t_c) = f(p)(t_d)$ so that
$t_c^{+} = t_d$ and $\fF_{f(p)}(t_c;t_d) \in \Bruhat_{(23);2}$.
By taking $\epsilon > 0$ sufficiently small,
we may assume that if 
$p \in s^{-1}(\Bruhat_{(123);6}) \smallsetminus V_2$
and  $t_0(p) -\epsilon < t_a < t_a^{+} < t_b < t_1(p) + \epsilon$
then $\fF_{f(p)}(t_a;t_b) \in \Bruhat_{(13);4}$.
Similarly, 
we may assume that if 
$p \in s^{-1}(\Bruhat_{(132);0}) \smallsetminus V_2$
and  $t_0(p) -\epsilon < t_a < t_a^{+} < t_b < t_1(p) + \epsilon$
then $\fF_{f(p)}(t_a;t_b) \in \Bruhat_{(13);1}$.

\begin{figure}[ht]
\begin{center}
\psfrag{ta}{$t_a$}
\psfrag{h+(ta)}{$h_{+}(t_a)$}
\psfrag{tb}{$t_b$}
\psfrag{h-(tb)}{$h_{-}(t_b)$}
\psfrag{tc}{$t_c$}
\psfrag{td}{$t_d$}
\epsfig{height=50mm,file=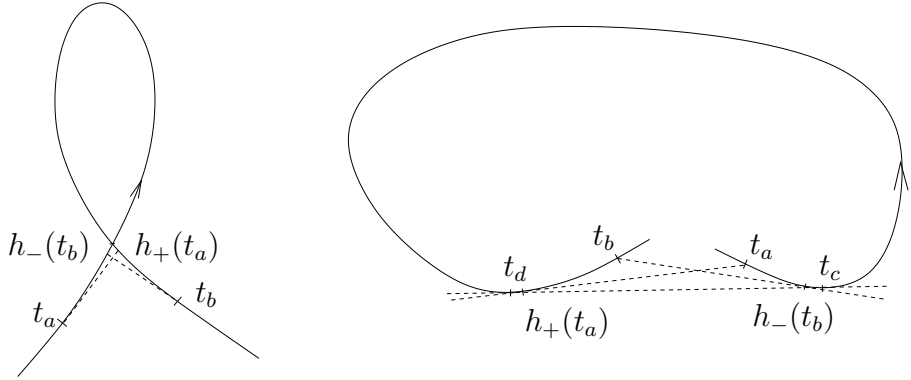}
\end{center}
\caption{The functions $h_+$ and $h_-$}
\label{fig:hh}
\end{figure}

If $p \in V_2$ and $t \ge t_c$ set $h_{+}(t) = t^{+}$;
if $t < t_c$,
let $\tilde t = h_+(t) \in (t_d,t_1+\epsilon)$ be uniquely defined
by $\fF_{f(p)}(t;\tilde t) \in \Bruhat_{(123);3}$;
more geometrically, draw a tangent geodesic to $\gamma$ at $t$:
the curve $\gamma$ intersects the geodesic at $\tilde t$.
The function $h_+: [t_c - \delta, t_c + \delta] \to [t_d, t_d + \tilde\delta]$
is continuous,
decreasing for $t < t_c$ and increasing for $t > t_c$
and satisfies $h_+(t_c) = t_d$ and $h_+'(t_c) = 0$.
Similarly, if $t \le t_d$ let $h_{-}(t) = t^-$;
otherwise, if $t > t_d$,
let $\tilde t = h_-(t) \in (t_0-\epsilon,t_c)$ be defined by
$\fF_{f(p)}(\tilde t; t) \in \Bruhat_{(132);6}$;
again,
the function $h_-: [t_d - \delta, t_d + \delta] \to [t_c - \tilde\delta, t_c]$
is continuous,
increasing for $t < t_d$ and decreasing for $t > t_d$
and satisfies $h_-(t_d) = t_c$ and $h_-'(t_d) = 0$
(see Figure \ref{fig:hh}).
By taking $\epsilon$ small we may assume that the functions
$h_+$ and $h_-$ are always $(1/2)$-Lipschitz,
i.e., that
$|h_+(t_a) - h_+(\tilde t_a)| \le |t_a - \tilde t_a|/2$ and
$|h_-(t_b) - h_-(\tilde t_b)| \le |t_b - \tilde t_b|/2$.

If $p \in V_2$,
$t_0 - \epsilon < t_a < t_a^+ < t_b < t_1 + \epsilon$
and $t_a > h_{-}(t_b)$ then $\fF_\gamma(t_a;t_b) \in \Bruhat_{(13),4}$.
Similarly,
if $t_b < h_{+}(t_a)$ then $\fF_\gamma(t_a;t_b) \in \Bruhat_{(13),1}$.
Thus, if
$\fF_\gamma(t_a;t_b) \notin \Bruhat_{(13),1} \cup \Bruhat_{(13),4}$
then
$t_a \le h_{-}(t_b) \le t_c < t_d \le h_{+}(t_a) \le t_b$.
Conversely,
if $t_a < h_{-}(t_b) < t_c < t_d < h_{+}(t_a) < t_b$
then $\fF_\gamma(t_a;t_b) \in \Bruhat_{(13),7}$.
But the Lipschitz condition means that if
$t_0 - \epsilon < t_a - \epsilon_b < t_a < t_a^+
< t_b < t_b + \epsilon_b < t_1 + \epsilon$
and $\fF_\gamma(t_a;t_b) \notin \Bruhat_{(13),1} \cup \Bruhat_{(13),4}$
then $\fF_\gamma(t_a - \epsilon_b;t_b + \epsilon_b) \in \Bruhat_{(13),7}$.

A similar construction holds in $U_4$,
only the Bruhat classes change and the geometric construction
of $h_{\pm}$ is different (see Figure \ref{fig:hh}):
to define $h_+(t_a)$ we construct a tangent line to $\gamma$
which passes through $t_a$.
Alternatively, we apply Arnold duality (\cite{Arnold}, \cite{SaSha})
to go from the $U_4$ to the $U_2$ scenario.
\qed

\section{Complicated curves}

\label{sect:Y}

Recall that a curve $\gamma \in \cL_I$
is \textit{complicated} if it belongs
to
\[ \cY_{Q} = \cL_{Q} \smallsetminus
\bigcup_{k} \cM_k. \]

\begin{lemma}
\label{lemma:looseY}
Let $K$ be a compact manifold and let $f: K \to \cL_Q$
be a continuous map.
If the image of $f$ is contained in $\cY_Q$
then $f$ is loose.
\end{lemma}

{\nobf Proof:}
Let $\tilde U_0 \subseteq K$ be the open set of elements
$p \in K$ for which $[0,0^{+}]$ is a good step.
More generally, let $\tilde U_j \subseteq K$ be the open set
of elements $p \in K$ for which $\Next_{f(p)}^{j+1}(0) < 1$
and $[\Next_{f(p)}^j(0), \Next_{f(p)}^{j+1}(0)]$ is a good step.
Since all curves $f(p)$ are complicated,
the sets $\tilde U_j$ cover $K$.
By compactness of $K$, there exists an integer $J$
such that $K = \bigcup_{j < J} \tilde U_j$.
Define compact sets $K_j \subset \tilde U_j$ such that
the interiors $U_j$ of $K_j$ already cover $K$.

Again by compactness, there exists $\epsilon_{J-1} > 0$
such that if $p \in K_{J-1}$ and
$|t_{J-1} - \Next_{f(p)}^{J-1}(0)| < \epsilon_{J-1}$
then $\Next_{f(p)}(t_{J-1}) < 1$ and $[t_{J-1}, \Next_{f(p)}(t_{J-1})]$
is a good step (for $f(p)$).
Define the continuous function $t_{J-1}: K_{J-1} \to [0,1]$
by $t_{J-1}(p) = \epsilon_{J-1}/2 + \Next_{f(p)}^{J-1}(0)$
Similarly, there exists $\epsilon_{J-2} > 0$
such that if $p \in K_{J-2}$ and
$|t_{J-2} - \Next_{f(p)}^{J-2}(0)| < \epsilon_{J-2}$
then $\Next_{f(p)}(t_{J-2}) < 1$,
$\Next_{f(p)}(t_{J-2}) < t_{J-1}(p)$ if $p \in K_{J-1}$
and $[t_{J-2}, \Next_{f(p)}(t_{J-2})]$ is a good step.
Proceed in this manner so that we have continuous
functions $t_j: K_j \to (0,1)$, $0 \le j < J$, such that
if $p \in K_j$ then $[t_j(p), \Next_{f(p)}(t_j(p))] \subset (0,1)$
is a good step for $f(p)$ and
if $p \in K_{j_1} \cap K_{j_2}$, $j_1 < j_2$,
then $\Next_{f(p)}(t_{j_1}(p)) < t_{j_2}(p)$.
Again by compactness, there exists $\epsilon > 0$ such that
$t_j(p) > \epsilon$, $\Next_{f(p)} < 1-\epsilon$ and
if $p \in K_{j_1} \cap K_{j_2}$, $j_1 < j_2$,
then $\Next_{f(p)}(t_{j_1}(p)) + \epsilon < t_{j_2}(p)$.

Thus if $p \in K_j$ then the interior of the interval
$\hat I_j = [ t_j(p) - \epsilon/3, \Next_{f(p)}(t_j(p)) + \epsilon/3 ]$
contains the good step $[t_j(p), \Next_{f(p)}(t_j(p))]$.
Notice furthermore that for given $p$
the intervals $\hat I_j$ are disjoint.
From Lemma \ref{lemma:goodarc} we can now 
define functions $a_j < b_j < c_j < d_j: K \to (0,1)$
such that, for all $p \in K$,
$a_j(p), b_j(p), c_j(p), d_j(p) \in \hat I_j(p)$
and, for $p \in K_j$,
$[b_j(p), c_j(p)] \subset [a_j(p), d_j(p)]$
is a good pair of arcs.
Define $I_j(p) = [b_j(p), c_j(p)]$
and $\tilde I_j(p) = [a_j(p), d_j(p)]$.

The strategy now is to deform curves in each interval $I_j$ independently.
More precisely, define
\begin{align*}
U_{j,1} &= \{ p \in K_j ;
\fF_{f(p)}(b_j(p); c_j(p)) \in \Bruhat_{(13),1} \}, \\
U_{j,4} &= \{ p \in K_j ;
\fF_{f(p)}(b_j(p); c_j(p)) \in \Bruhat_{(13),4} \}, \\
U_{j,7} &= \{ p \in K_j ;
\fF_{f(p)}(a_j(p); d_j(p)) \in \Bruhat_{(13),7} \}.
\end{align*}
By definition of good pair of arcs,
$U_j = U_{j,1} \cup U_{j,4} \cup U_{j,7}$.
Consider open sets $W_{j,\ell}$ and compact sets
$K_{j,\ell}$ such that $W_{j,\ell} \Subset K_{j,\ell} \Subset U_{j,\ell}$
and such that the sets $W_{j,\ell}$ cover $K$.

Let $f_0 = f$;
apply Lemma \ref{lemma:B147}
(with $\ell = 1$, $K_0 = K_{0,1}$, $W_0 = W_{0,1}$,
$t_0 = b_0$ and $t_1 = c_0$)
to define a homotopy from $f_0$ to another function $f_{0,1}$
with $f_{0,1}(p) = f_0(p)^{[b_0\#2;c_0\#2]}$ for all $p \in W_{0,1}$.
Apply the same lemma
(now with $\ell = 4$, $K_0 = K_{0,4}$, $W_0 = W_{0,4}$,
$t_0 = b_1$ and $t_1 = c_1$)
to define a homotopy from $f_{0,1}$ to $f_{0,4}$
with $f_{0,4}(p) = f_{0,1}(p)^{[b_0\#2;c_0\#2]}$ for all $p \in W_{0,4}$;
notice that since $K_{0,1}$ and $K_{0,4}$ are disjoint
the two constructions do not interfere with one another.
Apply Lemma \ref{lemma:B147} yet another time
(now with $\ell = 7$, $K_0 = K_{0,7}$, $W_0 = W_{0,7}$,
$t_0 = a_1$ and $t_1 = d_1$)
to define a homotopy from $f_{0,4}$ to $f_{0,7}$
with $f_{0,7}(p) = f_{0,4}(p)^{[a_0\#2;d_0\#2]}$ for all $p \in W_{0,7}$;
notice that $\fF_{f_{0,4}(p)}(a_0(p)) = \fF_{f(p)}(a_0(p))$
and $\fF_{f_{0,4}(p)}(d_0(p)) = \fF_{f(p)}(d_0(p))$.
Even though the interval $(a_0,d_0)$ contains the points $b_0$ and $c_0$,
the loops created in the two first steps were not destroyed
but merely pushed around by a projective transformation.
We may therefore define a homotopy from $f_{0,7}$ to $f_1$
such that
{if $p \in W_{0,1} \cup W_{0,4}$ then
$f_1(p) = \gamma^{[b_0\#2;c_0\#2]}$ for some $\gamma \in \cL_Q$;}
{if $p \in W_{0,7}$ then
$f_1(p) = \gamma^{[a_0\#2;d_0\#2]}$ for some $\gamma \in \cL_Q$.}

Repeat the process to define a homotopy from $f_1$ to $f_2$
and so on until $f_J$ such that, finally, for $f_J$ we have:
\begin{itemize}
\item{if $p \in W_{j,1} \cup W_{j,4}$ then
$f_J(p) = \gamma^{[b_j\#2;c_j\#2]}$ for some $\gamma \in \cL_Q$;}
\item{if $p \in W_{j,7}$ then
$f_J(p) = \gamma^{[a_j\#2;d_j\#2]}$ for some $\gamma \in \cL_Q$.}
\end{itemize}
Since the sets $W_{j,\ell}$ cover $K$
we have from Lemma \ref{lemma:looppox} that $f_J$ is loose.
\qed

We now prove that the inclusions $\cY_z \subset \cI_z$
are weak homotopy equivalences.

{\nobf Proof of Proposition \ref{prop:YI}: }
Let $K$ be a compact manifold of dimension $n$
and let $f: K \to \cI_z$ a continuous map.
From Lemma \ref{lemma:spread}, for sufficiently large $m$,
$f$ is homotopic to $f_1 = F_{2m} \circ f$
and the image of $f_1$ is contained in $\cL_z$.
Furthermore, the total curvature of $f_1(p) = (F_{2m} \circ f)(p)$
tends to infinity when $m$ tends to infinity.
We may therefore choose $m$ large enough so that
$\tot(f_1(p)) > 8(n+1)\pi$ for all $p \in K$.
In particular, the image of $f_1$ is disjoint from $\cM_{k}$
for all $k \le n+1$.
For $k > n+1$, the codimension of $\cM_{k}$ equals $2k-2 > 2n \ge n$:
by transversality, we may perturb $f_1$ to define a homotopic map
$f_2: K \to \cL_z$ whose image is disjoint from all submanifolds
$\cM_{k}$, so that the image of $f_2$ is contained in $\cY_z$.
In particular the maps $i_\ast = \pi_n(i): \pi_n(\cY_z) \to \pi_n(\cI_z)$
are surjective.

Conversely, let $B$ be a compact manifold of dimension $n+1$
with boundary $K = \partial B$.
Let $g: B \to \cI_z$ be a continuous map
with the image of $f = g|_{K}$ contained in $\cY_z$.
We prove that there exists a map $\tilde g: B \to \cY_z$
with $\tilde g|_{K} = f$.
Indeed, let $g_1 = F_{2m} \circ g: B \to \cI_z$.
Again from Lemma \ref{lemma:spread},
for sufficiently large $m$ we have that
the image of $g_1$ is contained in $\cL_z$.
From Lemma \ref{lemma:looseY},
$f: K \to \cL_z$ is homotopic (in $\cL_z$)
to $f_1 = g_1|_{K}$.
We therefore obtain a map
$g_2: B \to \cL_z$ with $g_2|_{K} = f$.
By Lemma \ref{lemma:Mkgohome},
there exists $g_3: B \to \cY_z \subset \cL_z$
with $g_3|_{K} = f$.
In particular the maps $i_\ast = \pi_n(i): \pi_n(\cY_z) \to \pi_n(\cI_z)$
are injective.
\qed

\section{Proof of Theorems \ref{theo:main} and \ref{theo:mainplus} \\
and final remarks}

\label{sect:theo}

We now have all the tools required to complete the proof
of Theorem \ref{theo:mainplus};
Theorem \ref{theo:main} is then a special case.

{\nobf Proof of Theorem \ref{theo:mainplus}: }
Recall that $\cM_k$ is either empty or
a contractible submanifold of codimension $2k-2$.
As in the proof of Lemma \ref{lemma:Mkgohome},
for each $k$, if $\cM_k \ne \emptyset$
let $V_k \subset \cL_z$ be
a closed tubular neighborhood of $\cM_k$
with interior $U_k \subset V_k$.
If $\cM_k = \emptyset$, let $U_k = V_k = \emptyset$;
also, for $k = 1$ if $\cM_k \ne \emptyset$ then
$\cM_k$ is a contractible connected component
and $U_k = V_k = \cM_k$.
As before, assume the sets $V_k$ to be disjoint.
If $\cM_k \ne \emptyset$,
let $\psi_k = (\psi_{k,1}, \psi_{k,2}): V_k \to \cH \times \DD^{2k-2}$
be a homeomorphism ($\cH$ is the Hilbert space).
Again, assume that $\cM_k = \psi_{k,2}^{-1}(\{0\})$.
Assume furthermore that $D_k = h_{2k-2}^{-1}(V_k) \subset \Ss^{2k-2}$ 
is a closed disk with smooth boundary.
Assume also that
$\psi_{k,1} \circ h_{2k-2}$ is constant equal to $0$ in $D_k$
and that $h_{2k-2}$ is a homeomorphism from $D_k$
to $\cD_k = h_{2k-2}(D_k) = \psi_{k,1}^{-1}(\{0\}) \subset V_k$. 
Let $\tilde\cY_z \subset \cY_z$ and $\tilde\cL_z \subset \cL_z$ be defined by
\[ \tilde\cY_z = \cL_z \smallsetminus \bigcup_{\cM_k \ne \emptyset} U_k, \quad
\tilde\cL_z = \tilde\cY_z \cup \bigcup_{\cM_k \ne \emptyset} \cD_k. \]
Since $\cH \times \Ss^{2k-3} \subset
\cH \times (\DD^{2k-2} \smallsetminus \{0\})$ is a deformation retract
it follows that so is $\tilde\cY_z \subset \cY_z$.
It thus follows from Proposition \ref{prop:YI}
that $\tilde\cY_z \approx \cI_z \approx \Omega\Ss^3$.
Similarly, since
$(\cH \times \Ss^{2k-3}) \cup (\{0\} \times \DD^{2k-2}) \subset
\cH \times \DD^{2k-2}$ is a deformation retract,
so is $\tilde\cL_z \subset \cL_z$.

If neither $z$ nor $-z$ is convex,
$\cM_k = \emptyset$ for all $k$ and
$\tilde\cY_z = \tilde\cL_z$ and we are done.
If $z$ (resp. $-z$) is convex then
$\cM_k \ne \emptyset$ for $k$ odd (resp. even).
Thus, for $z$ convex,
$\tilde\cL_z$ is obtained from $\tilde\cY_z \approx \Omega\Ss^3$
by gluing disks $\DD^0$ (i.e., adding a contractible connected component),
$\DD^4$, $\DD^8$, \dots 
Similarly, for $-z$ convex,
$\tilde\cL_z$ is obtained from $\tilde\cY_z \approx \Omega\Ss^3$
by gluing disks $\DD^2$, $\DD^6$, $\DD^{10}$, \dots 
The maps $h_{2k-2}$ guarantee that the spheres
along which these disks are being glued
are nullhomotopic in $\tilde\cY_z$
and therefore gluing a disk is (homotopically) equivalent
to gluing a sphere: the theorem follows.
\qed

The question of how the spaces $\cL_z$ fit together
still requires some clarification.
It should be noted that the map from $\cL$ to $\Ss^3$ taking
$\gamma$ to $\tilde\fF_\gamma(1)$ 
does not satisfy the homotopy lifting property (see also \cite{S1}).
In particular, we would like to gain a better understanding
of periodic solutions of a linear ODE of order $3$.

Finally, similar questions can be asked about curves in $\Ss^n$, $n > 2$
($\gamma$ is locally convex if $\det(\gamma(t), \ldots, \gamma^{(n)}(t)) > 0$);
in \cite{SaSha} we show a few results about these spaces.

It would also be interesting to investigate the homotopy type
of spaces of curves with bounded geodesic curvature.
In \cite{SZ} and \cite{Z} some first results are proved.

\bigbreak


\bigskip\bigskip\bigbreak

{

\parindent=0pt
\parskip=0pt
\obeylines

Nicolau C. Saldanha, PUC-Rio
saldanha@puc-rio.br; http://www.mat.puc-rio.br/$\sim$nicolau/



\smallskip

Departamento de Matem\'atica, PUC-Rio
R. Marqu\^es de S. Vicente 225, Rio de Janeiro, RJ 22453-900, Brazil

}

\end{document}